\documentclass{gtart}
\usepackage{amssymb}
\newcommand{\inv}{^{-1}}              
\newcommand{\basl}{\backslash}
\newcommand{\N}{{\mathbb N}}    
\newcommand{\Z}{{\mathbb Z}}
\newcommand{\R}{{\mathbb R}}
\def\kp{k\kern-.07em/\kern-.07em p\kern.16em }
\def\mlap#1{\hbox to0pt{\hss#1\hss}}%
\def\k#1{\kern#1truemm}%
\font\csc=cmcsc10
\let\ninepoint=\small
\def\unit{30 }%
\let\df=\sl 
\def\mlap#1{\hbox to0pt{\hss#1\hss}}%
\def\lwd#1{\hfilneg\rlap{\hbox{#1}}\hfil}
\def\frac#1#2{\hbox{$#1\over#2$}}%
\def\udlap#1{\setbox0\hbox{#1}%
  \vbox to0pt{\vss\copy0\kern-\dp0}\setbox0\null}%
\def\nf{\nu}%
\def\yconj#1{w_{#1}^{-1}y_{#1}^{\epsilon_{#1}}w_{#1}}%
\def\pq{\hbox{$*_p\Z_q$}}%
\def\pretend{\rlap{\phantom{$\Big($}}}%
\def\ei{{\epsilon_i}}%
\def\rnbs{\kern-.06em}%
\input epsf

\def\rllist{}%
\def\Mon{%
  \def\rllist{}%
}%
\def\Moff{%
  \rllist
  \def\rllist{}%
}%
\def\Mrl#1#2{%
  \special{ps:/a {} def}%
  \rlap{#2}%
}%

\def\diagplease#1{\special{" MTGdict begin startdiag #1 enddiag end}}%
\def\diagbox#1{\vbox to.7truein{\vfill\hbox to.64truein{\hfill
  \rlap{\diagplease{#1}}\hfill}\vfill}}%
\def\mdiag#1{\matrix{\hbox{\diagbox{#1}}}}%

\newtheorem{lemma}{Lemma}[section]
\newtheorem{theorem}[lemma]{Theorem}
\newtheorem{prop}[lemma]{Proposition}
\newtheorem{cor}[lemma]{Corollary}
\newtheorem{definition}[lemma]{Definition}
\newtheorem{lem}[lemma]{Lemma}

\begin{document}
%
%
\input gtoutput
\volumenumber{2}\papernumber{3}\volumeyear{1998}
\pagenumbers{31}{64}\published{21 March 1998}
\proposed{Cameron Gordon}\seconded{Joan Birman, Walter Neumann}
\received{4 August 1997}
\accepted{19 March 1998}

\title{A natural framing of knots} 

\author{Michael T Greene\\Bert Wiest}

\address{Mathematics Institute\\University of Warwick\\Coventry CV4 7AL, UK\\
\smallskip\\\rm Email:\stdspace\tt mtg@maths.warwick.ac.uk\stdspace{\rm 
or}\stdspace mtg@uk.radan.com\\bertold@maths.warwick.ac.uk}
\asciiaddress{Mathematics Institute, University of Warwick, Coventry CV4 7AL, 
UK. 
Email: mtg@maths.warwick.ac.uk, mtg@uk.radan.com, bertold@maths.warwick.ac.uk}

\begin{abstract}
Given a knot $K$ in the 3--sphere, consider a singular disk bounded by
$K$ and the intersections of $K$ with the interior of the disk. The
absolute number of intersections, minimised over all choices of
singular disk with a given algebraic number of intersections, defines
the {\sl framing function} of the knot.  We show that the framing
function is symmetric except at a finite number of points.  The
symmetry axis is a new knot invariant, called the {\sl natural
framing} of the knot.  We calculate the natural framing of torus knots
and some other knots, and discuss some of its properties and its
relations to the signature and other well-known knot invariants.
\end{abstract}
\asciiabstract{%
Given a knot K in the 3-sphere, consider a singular disk bounded by
K and the intersections of K with the interior of the disk. The
absolute number of intersections, minimised over all choices of
singular disk with a given algebraic number of intersections, defines
the framing function of the knot.  We show that the framing
function is symmetric except at a finite number of points.  The
symmetry axis is a new knot invariant, called the natural
framing of the knot.  We calculate the natural framing of torus knots
and some other knots, and discuss some of its properties and its
relations to the signature and other well-known knot invariants.}



\maketitlepage

Let $K\co S^1 \to S^3$ be an unoriented knot. Let $D$ be the 2--disk.
We define a {\sl compressing disk of} $K$ to be a map $f\co D\to S^3$
such that $f|_{\partial D}=K$ and such that $f|_{{\rm int}(D)}$ is 
transverse to $K$. 
Then $f|_{{\rm int}(D)}$ has only finitely many intersections with the knot. 
We call the intersection points the {\sl holes} of the compressing disk, and
denote their number by $n(f)$. So $n(f)=|\{ f\inv(K) \cap {\rm int}(D)\}|$.
One rather crude invariant of the knot $K$ is the {\sl knottedness} 
$$L(K):= \min\{n(f)\ |\ f \mbox{ a compressing disk}\},$$
which was first considered in \cite{Pannwitz}. 
(The English term `knottedness' was taken from \cite{BuZi}.)

We note that holes can occur with two different signs (depending on
the direction in which $K$ pierces $f(D)$), so we can refine the above
invariant 
by defining the {\sl framing function} $n_K\co \Z \to \N$ as
follows: for a given $k\in \Z$ we minimize the absolute number of
holes among all compressing disks with algebraically $k$ holes.
We shall see that (except at finitely many points in $\Z$) this
function is symmetric around some value of $k$, and we call this
`asymptotic symmetry axis' $k$ the {\sl natural framing} of $K$.
The aim of this paper is to determine the natural framing of certain classes
of 
knot, and to study its properties and relations with other knot invariants.

In section 1 we define the natural framing of knots and show that it shares
many properties with the signature. In section 2 we define the natural framing
of each component of a link, and calculate it for a number of links.
It seems that for most ``simple'' knots the natural framing number is even;
in section 3, however, we show that knots with odd natural framing do exist.
In section 4 we prove that the natural framing of the $(p,q)$-torus knot
is $-(p-1)(q-1)$. In section 5 we summarize 
all the information we have about
the natural framing of knots with up to seven crossings.

\section{Definitions and general results}

Let $K\co S^1 \to S^3$ be a knot, and let $f\co D \to S^3$ be a
compressing disk. By a standard general position argument we can
assume that the singularities of a compressing disk $f$ are all
transverse self-intersectons which may end either in $K$ or in the
branch points of Whitney umbrellas (for details see
\cite{umbrellas}). In what follows, we will, without further
explanation, talk about the double lines or the Whitney umbrellas of a
compressing disk.  The following result is from \cite{Pannwitz}.
\begin{lem}\label{noWU} \sl Let $f\co D \to S^3$ be a compressing disk. 
Then there exists a compressing disk $f'\co D \to S^3$, where $f'$ is an
immersion, such that $n(f')\leqslant n(f)$.
\end{lem}
\begin{proof} We homotope $f$ into general position. Now we only have to get rid
of Whitney umbrellas. The double line starting at a Whitney umbrellas either
ends in another Whitney umbrella or in $K$. We move the branch point
along the double line, shrinking the double line until it doesn't contain
triple points. In the first case we then perform a 
surgery along the double line; this eliminates the double line and
the two Whitney umbrellas. In the second case we 
slide the knot over the branch point
of the Whitney umbrella; this leaves the knot type unchanged, and
eliminates the  Whitney umbrella, the double line, and one hole. (Also, 
it changes the 
framing represented by the compressing disk by $\pm 1$.) 
Applying this process to every branch point gives the result. \end{proof}
%
\begin{prop}[Pannwitz]For any knot $K$, $L(K)$ is even. 
\end{prop}
\begin{proof} Let $f\co D\to S^3$ be a compressing disk with $n(f)=L(K)$.
As we have just seen, we may assume that $f$ is an immersion, so it has 
no Whitney umbrellas. Every hole of $f$ is the beginning of a double point
line. Now this line must end in another hole (not a Whitney umbrella); so
there must be an even number of holes.\end{proof}

{\bf Remark}\stdspace For any non-trivial knot $K$ we have $L(K)\geqslant 2$. This 
follows from Dehn's Lemma and the fact that $L(K)$ is even.
\begin{theorem}\label{Ladditive}The knottedness is additive under 
connected sum, ie
if $K_1$ and $K_2$ are knots then
$L(K_1\#K_2)=L(K_1)+L(K_2)$. 
\end{theorem}
\begin{proof} It is obvious that $L(K_1\#K_2) \leqslant L(K_1)+L(K_2)$; we need to
prove the opposite inequality.

Choose an embedded sphere $S$ in 
$S^3$ which intersects the knot in only two points, splitting 
$K_1$ from $K_2$. Let $f\co D \to S^3$ be a compressing disk of $K_1\#K_2$ 
with only $L(K_1\#K_2)$ holes. Make the map $f$ transverse 
to $S$ without changing the number of holes. $f^{-1}(S)$ consists of one arc 
$A$ in $D$ connecting 
two points on $\partial D$, and a number of disjoint circles in 
$D$. Choose an outermost one of those, and call it $c$.
Let $C$ be the disk bounded by $c$. Say $c$ represents an element 
$[k]\in H_1(S\basl(S \cap K)) \cong H_1(S^3\basl(K_1\#K_2))
\cong {\Z}$; then there are at least $k$ holes inside $C$.
$S\basl(S \cap K)$ is a sphere with two holes, so
we can replace $f|_C$ by a map whose image lies entirely in $S$, and
which has precisely $k$ intersections with $K$.
After pushing the image away from $S$ by a small homotopy, we have
a map $f'$ such that $f'|_C$ has precisely $k$ holes and maps no point 
of $C$, or a small neighbourhood of $C$, to $S$.

This construction has replaced the compressing disk $f$ by a
compressing disk $f'$ with at most as many holes (and the same
framing). Since by hypothesis $n(f)$ is minimal, we have $n(f')=n(f)$.
Applying the construction to all outermost circles of intersection
yields a compressing disk with $L(K_1\#K_2)$ holes, which has the
property that the arc $A$ is mapped to $S$, one of the two components
of $D\basl A$ is entirely mapped to one of the two components of
$S^3\backslash S$, and the other component of $D\basl A$ is mapped to the
other component of $S^3\backslash S$.  This gives rise to compressing
disks for the knots $K_1$ and $K_2$ with $n_1$ and $n_2$ holes
respectively such that $n_1+n_2=L(K_1\#K_2)$.\end{proof}

Let $d$ be the inner boundary of a small neighbourhood of $\partial D$ in $D$;
then $d$ is a curve in $D$ `close to $\partial D$'. $f$ maps $d$ to a 
longitude of the knot. If we orient $d$ and $\partial D$ in the 
same way, then $f(d)$ and $f(\partial D)=K$ have a linking number 
$[f(d)]\in H_1(S^3\basl K)\cong \Z$. We denote this number by $k(f)$.
So $d$ represents the {\sl framing} $k(f)$ of $K$.
Geometrically, this framing can be obtained as follows: choose an orientation 
of $D$; this determines an orientation of $\partial D$ and hence of $K$. 
For $i\in \{1,\ldots,n(f)\}$ let $x_i\in D$ be the $i$th hole of $f$. Let
$\sigma(i)=1$ if a positive basis of $T_{f(x_i)}f(D)$ followed by a 
positive tangent vector of $K$ at $f(x_i)$ forms a positive basis of
$T_{f(x_i)}S^3$, and let $\sigma(i)=-1$ otherwise.
Then $k(f)=\sum_{i=1}^{n(f)}\sigma(i)$.
\ppar

To every knot $K$ we can associate a function $n_K\co {\Z} \to
\N$ which we call the {\sl framing function} as follows: 
$$ n_K(k'):=\min\{n(f)\ |\ f \mbox{ a compressing disk with }
 k(f)=k'\} $$
Notice that $L(K)=\min n_K$. 
\begin{prop}\label{nproperties} \sl For any knot $K$, the function $n_K$ has
the following properties.
\begin{itemize}
\item[{\rm (i)}] $n_K(k) \geqslant |k|$ for all $k\in \Z.$
\item[{\rm (ii)}] $n_K$ maps even numbers to even numbers and odd numbers to 
odd numbers.
\item[{\rm (iii)}] `Continuity': For any $k\in \Z$ we have 
$n_K(k+1)=n_K(k) \pm 1$.
\item[{\rm (iv)}]  If $k\in \Z$ is odd, then $n_K(k)=
\min\{n_K(k-1),n_K(k+1)\}+1$. In particular, the function $n_K$ is 
completely determined by its values on even numbers.
\end{itemize} 
\end{prop}

\proof (i) is obvious.
\begin{itemize}
\item[(ii)] follows from the fact that $\sum_{i=1}^n \sigma(i)$ 
(where $\sigma(i)\in \{-1,1\}$) is even if and only if $n$ is even.
\item[(iii)] By (ii), $n_K(k+1)\ne n_K(k)$. So for definiteness say
$n_K(k+1) > n_K(k)$. Let $f$ be a compressing
map with $n_K(k)$ holes and framing $k$. By artificially introducing
a Whitney umbrella we can obtain a compressing disk with $n_K(k)+1$
holes and framing $k+1$.
\item[(iv)] Suppose $f$ is a compressing disk with odd framing number. 
Then $f$ has 
at least one Whitney umbrella. If we remove this as in Lemma 1.1, we obtain 
another compressing disk with framing number $k(f) \pm 1$ and $n(f)-1$ holes.
It follows that for every odd $z\in \Z$ we have either
$n_K(z+1)<n_K(z)$ or $n_K(z-1)<n_K(z)$. This implies (iv).\endproof
\end{itemize} 
\medskip
Let $B$ in $S^3$ be a double-point line of a compressing disk $f$.
Suppose that $B$ does not end in a Whitney umbrella, and suppose $B$
is not a closed curve.  Then $f^{-1}(B)$ consists of two lines in $D$
(not necessarily disjoint or embedded).  There are two possible cases:
either one of them connects two holes and the other connects two
points on $\partial D$, when $B$ is a {\sl ribbon singularity}; or each of
the two lines connects one hole with one point on $\partial D$, when
$B$ is a {\sl clasp singularity}. 
We call a clasp singularity {\sl positive} or {\sl negative} if it ends at
two positive or negative holes, respectively. 
We say a clasp or ribbon singularity is
{\sl short} if its double-point line meets no triple points.

\smallskip

{\bf Conjecture}\stdspace For any nontrivial knot $K$ we have $n_K(0)
\geqslant 4$.

\smallskip

This is a strengthened form of
Dehn's Lemma.  Dehn's Lemma states that the only knot in $S^3$ which has
a compressing disk without holes is the unknot.  Our conjecture 
asserts that this remains true under the weakened hypothesis that the
compressing disk has 2 intersections with $K$ of opposite sign (ie
one ribbon singularity $B$).  This is easy to show with the extra
hypothesis that $B$ is short.

\begin{definition}\sl The function $n_K$ gives rise to a framing of the knot:
$$
\nu(K):= \lim_{k\to \infty} {n_K(-k)-n_K(k) \over 2}
$$
is called the {\rm natural framing} of $K$.
Alternatively, we can define $\nu(K)$ to be the unique integer such that there
exists an $N\in \N$ with $n_K(\nu(K)-k)=n_K(\nu(K)+k)$ for all
$k \geqslant N$.\end{definition}
\begin{lem}\sl This is well-defined.\end{lem}
\begin{proof} For $k\in \N$ let $a_k=n_K(-k)-|{-}k|$ and let $b_k=n_K(k)-|k|$.
We have $(n_K(-k)-n_K(k))/2=(a_k-b_k)/ 2$. Moreover, by 
Proposition \ref{nproperties} (i) and (iii), both $a_k$ and $b_k$ are 
decreasing sequences in $2\N$. Therefore there exists an $N\in \N$ such that 
for $k\geqslant N$ the sequences $(a_k)$ and $(b_k)$ are constant. It follows 
that the limit exists. Furthermore, by \ref{nproperties} (ii) we have that
$n_K(-k)-n_K(k)$ is even for all $k\in \N$, so $\nu(K)$ is indeed an
integer.\end{proof}

The natural framing is an `asymptotic symmetry axis' of the framing function.
As a first example, we look at the framing function of the 
figure-of-eight knot $4_1$ (from the table in \cite{Rolfsen}).   
Only the value $n_{4_1}(0)=4$
\begin{figure}[htb]
$$
\Mon
\epsfbox{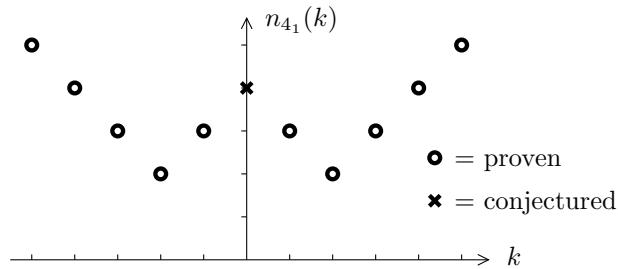}%
\Mrl{=conjectured}{\ninepoint\smash{\raise1pt\hbox{ = conjectured}}}%
\Mrl{=proven}{\ninepoint\smash{\raise1pt\hbox{ = proven}}}%
\Mrl n{\ninepoint\smash{\raise1pt\hbox{$n_{4_1}(k)$}}}%
\Mrl k{\ninepoint$k$}%
\Moff
$$
\vspace{-.8cm}
\caption{The function $n_{4_1}$}
\end{figure}
is conjectured. The fact that $n_{4_1}(2)=2$ follows from an easy 
construction (a compressing disk with one positive clasp 
and no other singularities). Similarly we see that $n_{4_1}(-2)=2$. 
The rest of the proof follows immediately from Proposition \ref{nproperties}.
We observe that $\nu(4_1)=0$. More generally we have:
\begin{prop}\label{nmirror} \sl Let $K$ be a knot, and let $mK$ be its 
mirror image. 
Then $n_K(k)=n_{mK}(-k)$ for all $k\in \Z$. In particular,
$\nu(K)=-\nu(mK)$, and if $K$ is amphicheiral then 
$\nu(K)=0$. \end{prop}
\begin{proof} Let $f$ be a compressing disk with framing $k$ and with $n_K(k)$ holes. 
Then $m\circ f$, where $m$ is the mirror map, is a compressing 
map of the knot $mK$ with framing $-k$ and also $n_K(k)$ holes.
It follows that $n_{mK}(-k)\leqslant n_K(k)$. The
opposite inequality is proved in the same way.\end{proof}

%
So the natural framing $\nu(K)$ changes sign under taking the mirror image, 
a property it shares with the well-known signature $\sigma(K)$ (see, for 
example, \cite{Rolfsen}). 
There is an even closer relation:  
\begin{prop}\sl If a knot $K$ has a compressing disk with positive clasps
and closed double-point lines, but no negative clasps or ribbons, then 
$\nu(K)\geqslant0$.  If all the clasps are short, then we 
have in addition $\sigma(K)\leqslant0$.\end{prop}
\begin{proof} If there exists a compressing disk with say $c$ positive clasps
and no negative clasps or ribbons then we have $n_K(k)=k$ for all
$k\geqslant 2c$.  It follows that $\nu(K)\geqslant0$.  If in addition
all these clasps are short, then we can unknot $K$ by $c$ negative
crossing changes. According to \cite{Giller} it follows that
$-2c\leqslant\sigma(K)\leqslant0$.\end{proof}

As an application, we can prove that for a large family of knots the
natural framing and the signature are both zero.  We only have to
construct a compressing disk with only short positive clasps and no
other singularities, and another compressing disk with only short
negative clasps.

\begin{figure}[htb]
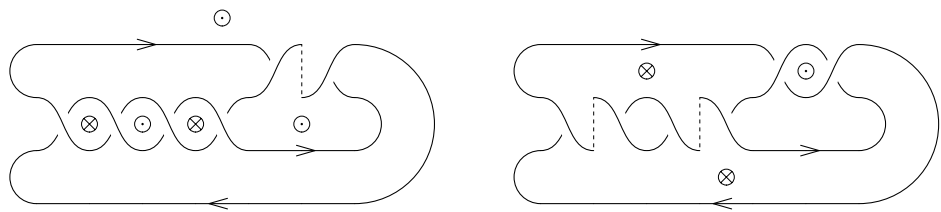

%
\vglue8truemm
\vbox to1truein{%
  \vfill
  $$
  \mdiag{%
gsave -90 rotate 
0 -5 tru
top sig1-nobl 0 0 up 1 0 90 arrw sig1-notr over1
 -4 ht -2 div up sig2 0 0 down 3 1 -90 arrw
 sig2 0 0 up -3 0 90 arrw sig2 0 0 down sig2 bot
grestore
gsave -90 rotate
0 14 tru
top sig1- -2 0 up 1 0 90 arrw sig1-
 2 ht -2 div down sig2nobr 3 1 -90 arrw
 sig2notl over2 -2 0 down -3 0 90 arrw sig2nobr sig2notl over2 bot
grestore
  }%
  $$
  \vfill
}%

\vspace{-.3cm} 

\caption{Two different compressing disks of the knot $6_1$} \label{mtg61}
\end{figure}

{\bf Examples}\stdspace(1)\stdspace If $K$ is one of the so-called `twist knots' $4_1,
6_1, 8_1, 10_1$ etc then $\nu(K)=\sigma(K)=0$. The case of the 
stevedore's knot $6_1$ is illustrated in Figure \ref{mtg61}, and
the other cases are similar. 
On the left we see a compressing disk with one short clasp singularity 
with positive sign (indicated by the dashed lines), proving that 
$n_{K_4}(2)\leqslant 2$. On the right we see a compressing 
disk with two short clasp singularities of negative sign, 
proving that $n_{K_4}(-4)\leqslant 4$.\medskip 

\begin{figure}[htb]
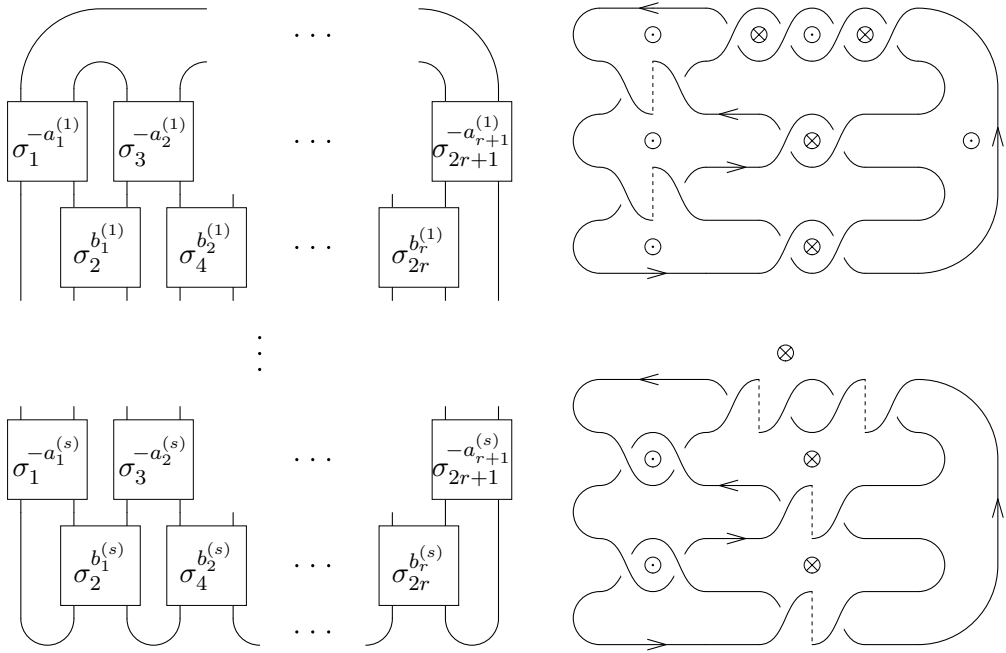

%
%
\bigskip
\vbox to3.3truein{\vfill
  \hbox to\hsize{\hfil\hfil\hfil\diagplease{toprightdiag}\k2\hfil}%
                \vfill
  \hbox to\hsize{\hfil
\Mon
\k2\diagplease{generaldiag}%
\Mrl1{\mlap{$\sigma_1^{\!-a_1^{(1)}}$}}%
\Mrl2{\mlap{$\sigma_3^{\!-a_2^{(1)}}$}}%
\Mrl3{\mlap{$\sigma_{2r+1}^{\!-a_{r+1}^{(1)}}$}}%
\Mrl4{\mlap{$\sigma_2^{b_1^{(1)}}$}}%
\Mrl5{\mlap{$\sigma_4^{b_2^{(1)}}$}}%
\Mrl6{\mlap{$\sigma_{2r}^{b_r^{(1)}}$}}%
\Mrl7{\mlap{$\sigma_1^{\!-a_1^{(s)}}$}}%
\Mrl8{\mlap{$\sigma_3^{\!-a_2^{(s)}}$}}%
\Mrl9{\mlap{$\sigma_{2r+1}^{\!-a_{r+1}^{(s)}}$}}%
\Mrl a{\mlap{$\sigma_2^{b_1^{(s)}}$}}%
\Mrl b{\mlap{$\sigma_4^{b_2^{(s)}}$}}%
\Mrl c{\mlap{$\sigma_{2r}^{b_r^{(s)}}$}}%
\Moff
                      \hfil\hfil\hfil}%
                \vfill
  \hbox to\hsize{\hfil\hfil\hfil\diagplease{botrightdiag}\k2\hfil}%
                \vfill
}%
\caption{A family of knots whose natural framing and signature are 
zero}\label{mtgplat}
\end{figure}

(2)\stdspace This construction can be generalized. Consider the family
of $r$--bridge knots indicated in Figure \ref{mtgplat}, 
where the $a_j^{(i)}$ and $b_j^{(i)}$ are all
non-negative and even.  On the right, choosing a particular
12--crossing knot as an example, we see two compressing disks; the
first has only short positive clasps, the second only short negative
clasps.


\medskip

{\bf Conjecture}\stdspace The natural framing of the twist knots $3_1$, $5_2$,
$7_2$, $9_2$, $\ldots$ is $2$.  More precisely, we conjecture that for
such knots $n(k)=2+|k-2|$.  The values for $k>0$ are easy to prove; the
values for $k\leqslant0$ seem hard.

Even more generally, we can consider the
${4ml+1\over 2l}$--two-bridge knot (in the notation of
\cite{BuZi}, Chapter 12).  This is the rational knot $C(2m \ 2l)$
in the notation of Conway \cite{Conway}, and for $l=1$ we get twist
knots.  We conjecture that for $m,l\in\Z^+$ the natural framing of
this knot is $\min(2m,2l)$.

\begin{lem}\label{sumle}\sl
Let $K_1$ and $K_2$ be knots. Then
$$
n_{K_1\#K_2}(k)=\min_{k'\in\Z}\,\bigl(n_{K_1}(k')+
n_{K_2}(k-k')\bigr).
$$
\end{lem}
\begin{proof}As in the proof of Theorem \ref{Ladditive} we see that for 
every compressing disk $f$ of $K_1\#K_2$ with framing $k$ and 
$n_{K_1\#K_2}(k)$ holes
we can find another compressing disk $f'$ with the same framing, the same
number of holes and only one intersection curve with a separating sphere. 
The result follows. \end{proof}
\begin{prop} \label{nyadditive} \sl
The natural framing is additive under connected sum; ie if $K_1$ and $K_2$
are knots then $\nu(K_1\#K_2)=\nu(K_1)+\nu(K_2)$.
\end{prop}
\begin{proof}
There exist $N,M \in \N$ and $c_1, c_2 \in \Z$ such
that $n_{K_1}(\nu(K_1)+k)=c_1+|k|$ and 
$n_{K_2}(\nu(K_2)+l)=c_2+|l|$ for all $k,l \in \Z$ with 
$|k|\geqslant N$ and $|l|\geqslant M$.

Now let $k\geqslant N$ and $l\geqslant M$. By Lemma \ref{sumle} we have
$$
\begin{array}{r@{}c@{}l}
n_{K_1\#K_2}(\nu(K_1)+\nu(K_2)+k+l)&{}\leqslant{}& 
   n_{K_1}(\nu(K_1)+k)+n_{K_2}(\nu(K_2)+l)\\
 &{}={}&c_1+k+c_2+l.
\end{array}
$$
Furthermore we have $n_{K_1}(\nu(K_1)+k+a)\geqslant c_1+k+a$ and
$n_{K_2}(\nu(K_2)+l-a)\geqslant c_2+l-a$ for all $a\in \Z$.
Lemma \ref{sumle} implies that 
$$
\begin{array}{r@{}c@{}l}
n_{K_1\#K_2}(\nu(K_1)+\nu(K_2)+k+l)
&{}={}&\min\limits_{a\in \Z}\>\Bigl(n_{K_1}(\nu(K_1)+k+a)\vspace{-2mm}\\
 &&\hspace{2cm} {}+n_{K_2}(\nu(K_2)+l-a)\Bigr)\\
 &{}\geqslant{}&c_1+k+c_2+l.
\end{array}
$$
We have proved that 
$n_{K_1\#K_2}(\nu(K_1)+\nu(K_2)+k)=c_1+c_2+|k|$ 
for all $k\in \Z$ with $k\geqslant N+M$. The case $k\leqslant -N-M$ is 
proved similarly.\end{proof}


\section{A natural framing of links}
In this section we define a natural framing for each component of a link.
By giving an example we prove that these framing numbers are {\it not} 
always even, and that they are {\it not} determined by the natural framings 
of the individual link components (regarded as knots) and their linking 
numbers. 

Let $L=L_1\cup\ldots\cup L_m\co S^1\cup\ldots\cup S^1 \to S^3$ be an 
unoriented link with $m$ components. Let $D$ be the 2--disk. We define a 
{\sl compressing disk of the $i$th link component}
$L_i$ ($i\in\{1,\ldots,m\}$) to be a map 
$f\co D\to S^3$ transverse to $L$ such such that $f|_{\partial D}=L_i$. 
Then $f|_{{\rm int}(D)}$ has only finitely many intersections with $L$. We 
call these intersection points the {\sl holes} of the compressing disk, and
denote their number by $n(f)$.  
We choose an orientation of $L_i$.  This induces an orientation of $D$.  We 
look at an intersection point of $f|_{{\rm int}(D)}$ with $L_i$.  We define 
such a hole to be {\sl positive} or {\sl negative}, depending on 
whether a positive basis of the tangent space to $D$ followed by a positive
tangent vector to $L_i$ forms a positive or a negative basis of $S^3$,
respectively.
This is well-defined (ie independent of the choice of orientation
of $L_i$).  We denote by $k(f)\in\Z$ the number of intersections of 
$f|_{{\rm int}(D)}$ with $L_i$ (not {\it all}\/ link components!), counted 
algebraically. Again, we can think of $k(f)$ as the framing of $L_i$ defined 
by $f$.

To every component $L_i$ of $L$ we can associate a function
$n_i\co {\Z} \to \N$ which we call the {\sl $i$th framing function} as
follows:
$$
n_i(k'):=\min\{n(f)\ |\ f\hbox{ a compressing disk of }L_i%
 \hbox{ with }k(f)=k'\}.
$$
Precisely as in the case of knots we define the {\sl natural framing of
the component $L_i$ of $L$} by
$$
\nu_i(L):=\lim_{k\to\infty}{n_i(-k)-n_i(k)\over2},
$$
and prove that this limit exists.
It is clear from the definition that the natural framing of each component
of $L$ is an integer multiple of ${1\over2}$. We claim:
\begin{lem}\sl The natural framing of each component of $L$ is an
integer.\end{lem}
\begin{proof} It suffices to show that for each $i\in\{1,\ldots,m\}$ we have
either $n_i(k) \equiv k$ (mod $2$) for all $k\in \Z$ or
$n_i(k)+1 \equiv k$ (mod $2$) for all $k\in \Z$. To see this,
let $f$ be a compressing disk with $n(f)=n(k')$, where $k'=k(f)$.
Equip all link components with an orientation, no matter which.
Let $\tilde{k}(f)$ be the number of 
intersections of $f|_{{\rm int}(D)}$ with {\it all} components of $L$,
counted algebraically. We have
$$ \tilde{k}(f) \equiv n(f) \mbox{ (mod } 2)$$
and
$$ k(f)=\tilde{k}(f)-\sum_{j\neq i} lk(L_i,L_j).$$
It follows that 
$$ n(f)-k(f) \equiv \sum_{j\neq i} lk(L_i,L_j) \mbox{ (mod } 2).$$
Since $\sum_{j\neq i} lk(L_i,L_j)$ is independent of $f$, the result
follows.\end{proof}

{\bf Examples}\stdspace (1)\stdspace  The trivial link on $m$ components.  Let $L$ be the
link consisting of $m$ unknotted, unlinked components.  It is easy to see 
that the framing function of each component is $n_i(k)=|k|$.  It follows 
that $\nu_i(L)=0$ for each $i\in \{1,\ldots,m\}$.
\medskip

(2)\stdspace The Hopf link.  The two link components $L_1$ and $L_2$ have linking 
number 1, so
any compressing disk of the first component has at least one intersection
with the second.  Therefore we have $n_1(k)\geqslant |k|+1$ for all $k\in \Z$.
It is easy to construct compressing disks of $L_1$ which have framing $k$
and $|k|+1$ holes, so $n_1(k)=|k|+1$ for all $k\in \Z$.  It follows that
$\nu_1(L)=0$, and similarly $\nu_2(L)=0$.
\medskip

\begin{figure}[htb]
$$
\Mon
\epsfbox{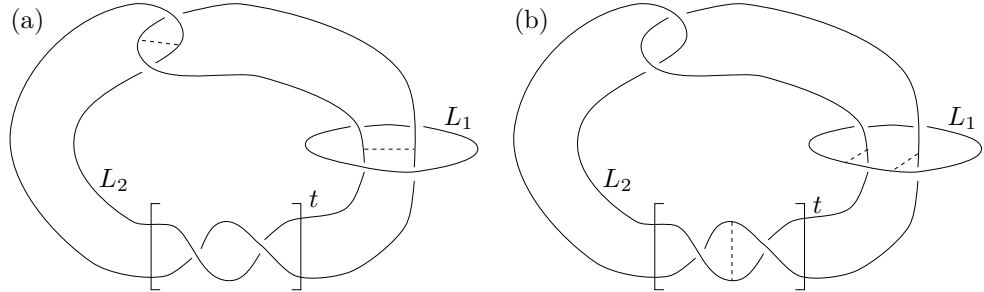}%
%
%
\Mrl{L1}{\ninepoint$L_1$}%
\Mrl{L1'}{\ninepoint$L_1$}%
\Mrl{L2}{\ninepoint$L_2$}%
\Mrl{L2'}{\ninepoint$L_2$}%
\Mrl{t}{\ninepoint$\!t$}%
\Mrl{t'}{\ninepoint$\!t$}%
\Mrl{a}{\ninepoint(a)}%
\Mrl{b}{\ninepoint(b)}%
\Moff
$$
\vspace{-.8cm}
\caption{Two different compressing disks of $L_2$}\label{oddlink}
\end{figure}

(3)\stdspace The link with two components shown in Figure \ref{oddlink} for any 
$t\geqslant 1$, consisting of an unknot and a twist-knot (eg for $t=2$ 
we get the knot $6_1$).
\begin{prop} \sl
We have $\nu_1(L)=0$ and $\nu_2(L)=1$
\end{prop}
\begin{proof} $L_1$, regarded only as a closed curve in $S^3\basl L_2$, represents a
nontrivial element of $\pi_1(S^3\basl L_2)$. Therefore any compressing disk
of $L_1$ has at least one intersection with $L_2$. Since the linking
number of $L_1$ and $L_2$ is 0, there must in fact be a minimum of 2 holes.
Therefore we have $n_1(k)\geqslant |k|+2$ for $k\in \Z$, and again equality
follows by construction. It follows that $\nu_1(L)=0$.

In order to visualize 
compressing disks of $L_2$, we draw their lines of self-intersections, and 
also the lines of intersection with the obvious compressing disk of $L_1$ 
which has no self-intersections and two intersection points with $L_2$.

There exists a compressing disk of $L_2$, indicated in Figure 
\ref{oddlink}(a), which is disjoint from $L_1$, and has one clasp singularity 
connecting two positive holes. This proves that $n_2(k)=k$ for $k\geqslant 2$.

On the other hand, there exists a compressing disk, indicated in Figure 
\ref{oddlink}(b), which has two intersections of opposite sign with $L_1$, 
and $t$ clasp singularities, each connecting two negative holes. This proves 
that $n_2(k)\leqslant |k|+2$ for $k\leqslant -2t$. 

Next let $f$ be a compressing disk with $k(f)\leqslant 0$. The image of $f$ 
has an even number of intersections with $L_1$, because the linking number 
of $L_1$ with $L_2$ is $0$. We distinguish two 
cases. If $f$ is disjoint from $L_1$, ie if the image of $f$ is contained
in the solid torus $S^3\basl L_1$, then
using a simple covering space argument it 
is easy to prove that $f$ has at least two positive holes.
So $n(f)\geqslant |k(f)|+4$ for all such compressing disks.
If, however, the image of $f$ has $2s$ intersections with $L_1$
($s\geqslant 1$) then $n(f)\geqslant |k(f)|+2s$. In either case, 
$n(f)\geqslant |k(f)|+2$. It follows that $n_2(k)\geqslant |k|+2$ for 
$k\leqslant 0$.

Altogether we have $n_2(k)=|k|+2$ for $k\leqslant -2t$, and therefore 
$\nu_2(L)=1$.\end{proof}

{\bf Example}\stdspace(4)\stdspace The Whitehead link, which is the case $t=0$ in Figure \ref{oddlink}. 
We have
$\nu_2(L)=1$, by precisely the same argument as in the case $t>0$.
The Whitehead link is isotopic to itself with the roles of $L_1$ and
$L_2$ interchanged.  Therefore $\nu_1(L)=1$.  (Note that the reasoning
behind the calculation of $\nu_1(L)$ for the case $t>0$ does not apply here,
since the {\it path} $L_1$ {\it is} contractible in $S^3\basl L_2$.)%
\ppar

These results are remarkable, because they show that $\nu_i(L)$, the
natural framing of the $i$th link component, is not determined
by the natural framing numbers of all individual link components
and their linking numbers. More precisely, $\nu_i(L)$ is not determined
by $\nu(L_1),\ldots,\nu(L_m)$ and $lk(L_r,L_s)$ ($r,s\in \{1,\ldots,m\}$).
%
Also, the natural framing numbers of link components can be odd. It is
not obvious that {\it knots} with odd natural framings exist, but
example (3) will lead to the construction of such a knot in the next
section.

\begin{prop}\label{linkcons}\sl
Let $L^{(1)}=L^{(1)}_1\cup\ldots\cup L^{(1)}_r$ and
$L^{(2)}=L^{(2)}_1\cup\ldots\cup L^{(2)}_s$ be links in $S^3$ with $r$
and $s$ components, respectively.  Let $L^{(3)}$ be the link obtained by
embedding $L^{(1)}$ and $L^{(2)}$ on either side of some embedded $S^2$ in 
$S^3$, and connecting $L^{(1)}_1$ and $L^{(2)}_1$ by a `band', as in the
construction of the connected sum of two knots.  $L^{(3)}$ has $r+s-1$
components, which we label such that $L^{(3)}_1$ is the one that
contains the `band'.  Then
$\nu_1(L^{(3)})=\nu_1(L^{(1)})+\nu_1(L^{(2)})$.
\end{prop}
\begin{proof} The proof is virtually identical to the proof of \ref{nyadditive}.
(Note that the link $L^{(3)}$ does not depend on the 
choice of the band.)\end{proof}

\begin{cor}\sl Let $L=L_1\cup\ldots\cup L_r$ be a link in $S^3$.
Then we can add one unknotted link component to $L$ in such a way that 
the natural framing of $L_1$ is increased by $1$, and such that the
natural framings of $L_2,\ldots,L_r$ remain unchanged.\end{cor}
\begin{proof} Both components of the Whitehead link are unknotted, so taking
the connected sum of $L_1$ with either of its components doesn't
change the type of $L$.  The result now follows from Example (4) and
Proposition \ref{linkcons}.\end{proof}


\section{A knot with odd natural framing}

For `simple' knots, eg knots with low crossing number, the natural framing
always appears to be an even number. In this section we exhibit a knot $K$
with $\nu(K)=1$. This knot is a satellite of a connected sum of three knots.
We do not know an atoroidal knot with odd natural framing number.

\begin{figure}[htb]
$$
\epsfbox{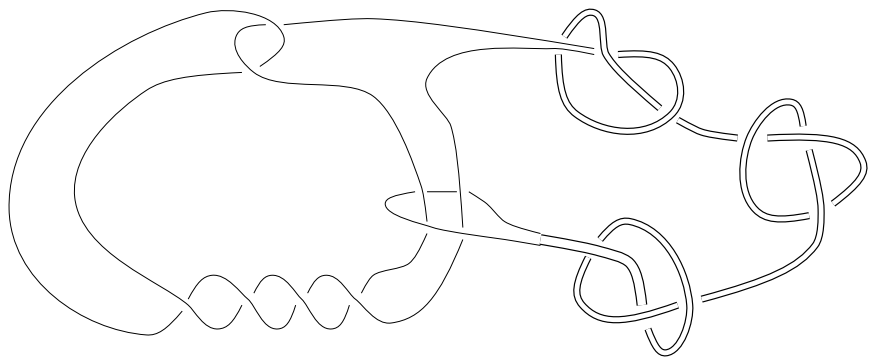}%
$$
\vspace{-.9cm}
\caption{The knot $K$ with $\nu(K)=1$}\label{oddkn}
\end{figure}
We need one more technical tool. Let $K$ be a knot, and let $f\co D \to S^3$
be any continuous map. We define a {\sl branched hole of} $f$ to be a point
$p\in D$ such that $f(p)\in K$ and such that $f$ maps the boundary of a 
small disk containing $p$ to some power of a meridian of $K$. So 
branched holes may just be transverse intersections of $f(D)$ with $K$,
but they may also be essentially nontransverse. 
If all intersections of $f(D)$ with $K$ are branched holes,
then we say $f$ is {\sl branched transverse}.

Let $\rho_n\co S^1 \to S^1$ be the standard map of degree $n$. Then any
branched transverse map $f\co D \to S^3$ with $f|_{\partial D}=K\circ \rho_n$
has at least one branched hole, by the loop theorem. More generally:
\begin{lem} \sl \label{branchedlem}
Let $K_1,\ldots,K_m$ be nontrivial knots. Then any branched transverse
map $f\co D \to S^3$ with $f|_{\partial D}=(K_1\#\ldots\# K_m)\circ \rho_n$
has at least $m$ branched holes.
\end{lem}
\begin{proof} The proof is similar to, but even simpler than, the proof of 
Theorem \ref{Ladditive}. For any given compressing disk $f$ there exists a
compressing disk $f'$ with no more branched holes than $f$ and with only
$m$ arcs (no closed curves) of intersection with a separating sphere.
Then $f'$ can be split into two disks, and the lemma follows inductively. \end{proof}

We are now ready to prove the main result of this section.
Let $K$ be the knot indicated in Figure \ref{oddkn}.
\begin{theorem}\label{mainth}The framing function of $K$ satisfies 
$n(k)=|k|+2$ for $k\geqslant 2$ and $n(k)=|k|+4$ for $k\leqslant -4$.
In particular, $\nu(K)=1$. \end{theorem}
\begin{proof}
It is easy to construct a compressing disk with $6$ negative and $2$ positive
holes, so $n(k)\leqslant |k|+4$ for $k\leqslant -4$ 
(see Figure \ref{optimcd}(a)). It is also easy to construct a compressing 
disk with $1$ negative and $3$ positive holes (see
Figure \ref{optimcd}(b)), 
\begin{figure}[htb]
$$
\epsfbox{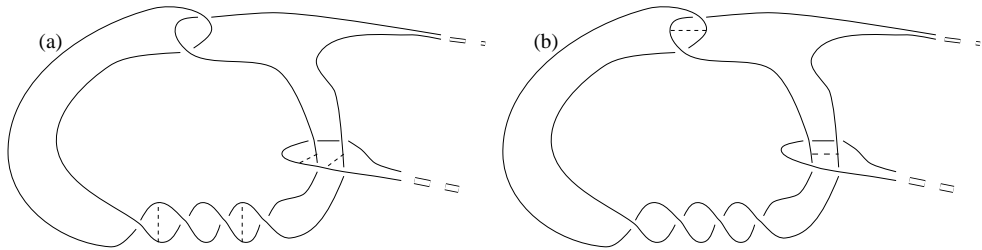}%
$$
\vspace{-.9cm}
\caption{The two `optimal' compressing disks of $K$}\label{optimcd}
\end{figure}
so $n(k)\leqslant |k|+2$ for $k\geqslant 2$. We want to prove that these
are in fact equalities, ie that every compressing disk of $K$ has
at least one negative and two positive holes.

Consider the knotted solid torus $S$ containing $K$, with a meridinal 
curve $c$ on its boundary, as indicated in Figure \ref{torusS}. 
The core of $S$ is a connected sum of three trefoil knots.
\begin{figure}[hbt]
$$
\epsfbox{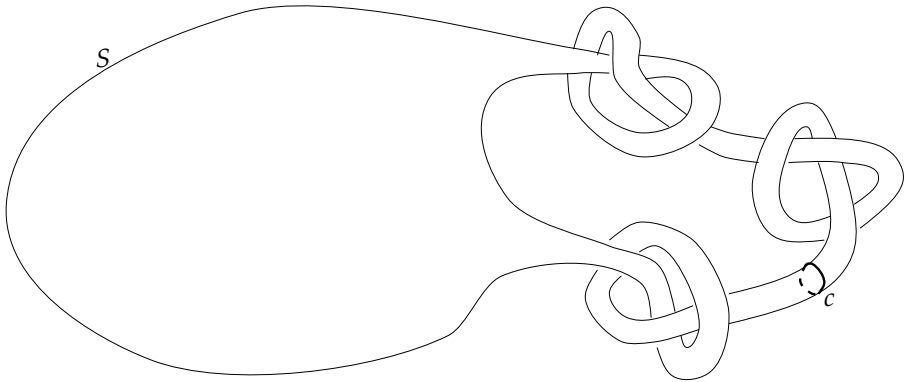}%
$$
\vspace{-.9cm}
\caption{The solid torus $S$ containing $K$}\label{torusS}
\end{figure}

Let $f\co D\to S^3$ be a compressing disk transverse to $\partial S$ with
$n(f)=n(k(f))$, ie $f$ has the minimal possible number of holes for
its framing. Two compressing disks $f_0$ and $f_1$ are called {\sl isotopic} 
if there is a homotopy $f_t$ ($t\in [0,1]$) which is fixed on $\partial D$
$$
f_t \co (D,D-f_0\inv(K)) \to (S^3,S^3-K).
$$
We can assume that among all disks isotopic to $f$, 
$f$ has the minimal number of intersections points $f(D)\cap c$.
\begin{lem} \sl \label{fcdisj} The compressing disk does not intersect 
the curve $c$.
\end{lem}
\begin{proof}[Proof of the Lemma]Assume it does. Then we look at intersection
lines $f(D)\cap \partial S$. Each such curve represents an element
$(m,l)\in H_1(\partial S) \cong \Z^2$, with $(1,0)$ corresponding to a 
standard meridian of $S$. We can assume that there are no 
inessential curves, ie no curves representing $(0,0)$, because we can
remove them by an isotopy of $f$. We call curves
representing $(m,0)$ ($m\in \Z\basl 0$) {\sl meridinal} curves, and all others
except the trivial one {\sl longitudinal} curves.
Since $|f(D)\cap c|$ is assumed minimal, meridinal curves are disjoint
from $c$, so there is at least one longitudinal curve.
The preimages of longitudinal curves are disjoint embedded circles in $D$,
and we let $\delta$ be an innermost one. Then $\delta$ bounds a disk $\Delta$
in $D$ such that $f|_\Delta$ has only meridinal intersections with 
$\partial S$. These meridinal curves are noncontractible in $S-K$, 
and they have linking number $0$ with $K$, so a disk in $\Delta$
bounded by the preimage of a meridinal curve contains at least one positive 
and one negative hole. By Lemma \ref{branchedlem} there are at least three
such meridinal curves. It follows that the framing disk has at least three
positive and three negative holes, contradicting the hypothesis that
the number of holes is minimal for its framing. \end{proof}

$S^3 \basl c$ is a solid torus; since, by Lemma \ref{fcdisj}, $f(D)$ is
disjoint from $c$, we can lift $K$ and the compressing disk $f$ to its 
universal cover. This is an open, infinite solid cylinder, and thus 
homeomorphic to $\R^3$ (see Figure \ref{oddKlift}). 

The preimage of $K$ under the covering space projection consist of a 
$\Z$--family of link components $\ldots,L_{-1},L_0,L_1,L_2,\ldots$, but 
we simply take away the link components 
$\ldots,L_{-3}$, $L_{-2},L_{-1}$ (see Figure \ref{oddKlift}). We denote by
$f'$ the lifting of the disk which sends $\partial D$ to $L_1$. $L_1$, 
regarded only as a closed curve, is noncontractible in $\R^3\basl L_2$ and 
has linking number zero with $L_2$, so $f'$ has at least one positive and 
one negative intersection with $L_2$. It follows that $f$ has at least one 
positive and one negative hole.

We can now also forget about the link components $L_i$ for $i\geqslant 2$,
and only consider $L_0$ and $L_1$.
\np                              
\phantom{\tiny .} \kern -\baselineskip \vspace{-3mm}
$$
\epsfbox{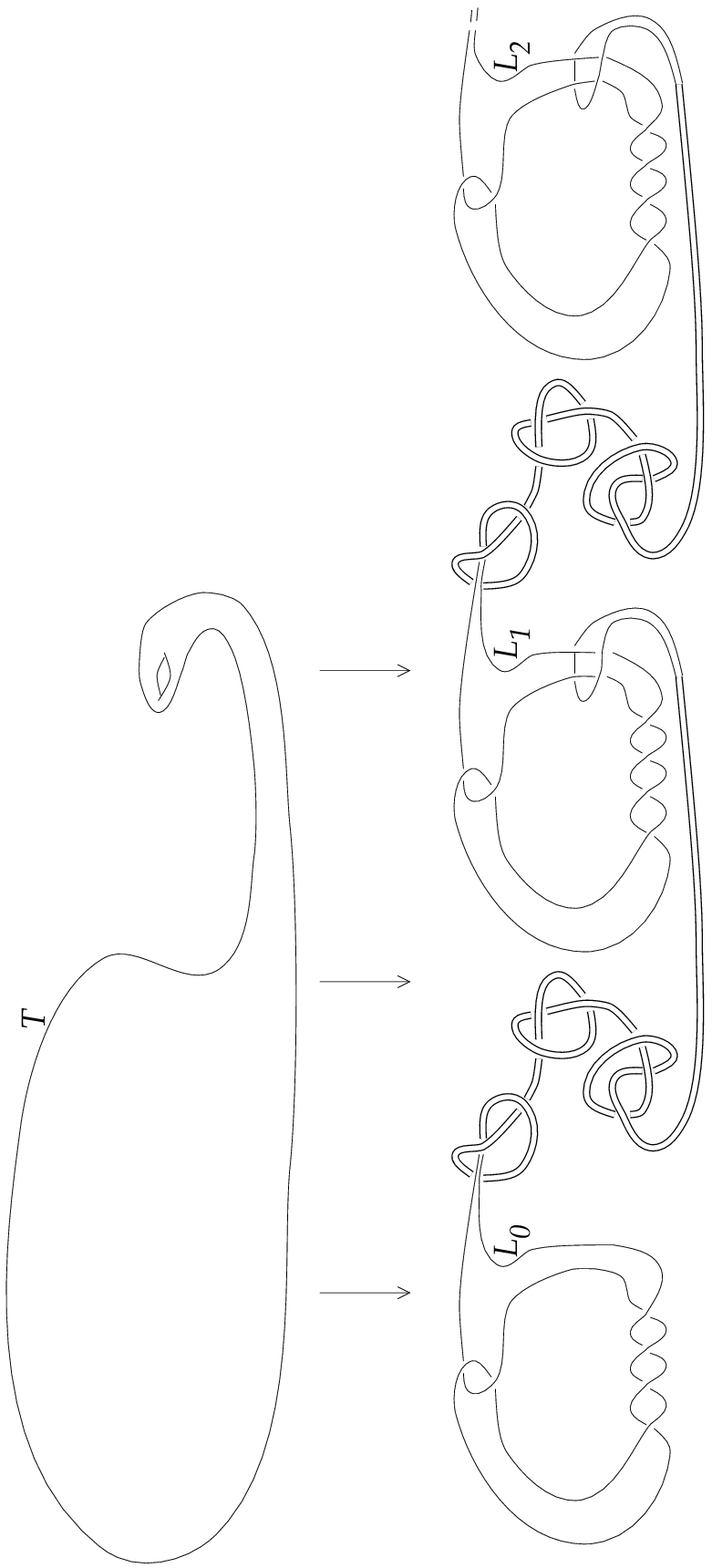}
$$
\vspace{-.9cm}
\begin{figure}[htb]
\caption{The universal cover of $S^3-c$}\label{oddKlift}
\end{figure}
\np
\begin{lem} \sl
The disk $f'$ has at least one positive intersection, either with $L_1$
or with $L_0$. 
\end{lem}
\begin{proof}[Proof of the Lemma]We embed a solid torus $T$ in $\R^3$ such that
it contains $L_0$, as indicated.
We can assume that $f'$ has the minimal number of intersection lines with 
$\partial T$ among all compressing disks isotopic to $f'$.
There are two possibilities to consider. Either $f'$ does not intersect $T$.
Then again by a covering space argument we see that $f'$ has at least two
positive intersections with $L_1$.
The other possibility is that $f'$ {\it has } intersections with 
$\partial T$. The preimages of these intersection lines are disjoint embedded 
circles in $D$, and we denote by $\gamma$ an innermost one of them; $\gamma$ 
bounds a disk $\Gamma\subseteq D$. Recall that $f'(\gamma)$ is assumed 
noncontractible in $\partial T$. So either $f'(\gamma)$ is a power of the 
$0$--longitude of $T$,
in which case $f'(\gamma)$ has at least one positive and one negative 
intersection with $L_1$; or $f'(\gamma)$ has nonzero linking number with 
the core of $T$, in which case the compressing disk must intersect $L_0$.
Since ${\rm lk}(L_1,L_0)=0$, we have at least one positive and one
negative intersection with $L_0$. \end{proof}

In summary, $f'$ has at least one negative and one positive intersection
with $L_2$, and at least one more positive intersection with $L_1$ or $L_0$, 
so $f$ has at least one negative and two positive holes. 
This completes the proof of Theorem \ref{mainth}.\end{proof}

The signature of a knot is always an even number.
So in particular, the absolute value of the natural framing and the signature 
of a knot can definitely be different.


\section{The natural framing of torus knots}

In this section we prove that the natural framing of $T(p,q)$, the 
$(p,q)$--torus knot, is $-(p-1)(q-1)$. More precisely, we exhibit a 
compressing disk of $T(p,q)$ with
$(p-1)(q-1)$ negative 
and no positive 
holes, and we prove that every
compressing disk must have at least $(p-1)(q-1)$ 
negative holes. A very different, more constructive proof of the special
case $q=2$ can be found in \cite{nat3}.

We start by reinterpreting the framing function.
Let $K\subseteq S^3$ be an oriented knot, and let
$G=\pi_1(S^3\basl K)$ be the knot group. A {\sl positive} or {\sl negative 
Wirtinger generator} of $G$ is a path from the basepoint to the boundary of a 
tubular neighbourhood of $K$, then once around a meridian of the tubular 
neighbourhood according to the right or left hand rule respectively,
and back along the first 
segment of path in the opposite direction.
Let $x_1,\ldots,x_m$ be positive Wirtinger generators
which together generate the
knot group.  Fix a path $\gamma$ from the
boundary of a tubular neighbourhood of the knot to the basepoint.
For $k\in \Z$ let $l_k\in G$ be the element represented by the 
path $\gamma\inv$ followed by the longitude with linking number $k$ with the 
knot followed by the path $\gamma$. 
\begin{lem}\sl Take $\sigma_+,\sigma_-\in\N$.  Let 
$k=\sigma_+-\sigma_-$ and $n=\sigma_++\sigma_-$.  Then the following 
statements are equivalent:
\begin{enumerate}
\item[{\rm (i)}] There exists a compressing disk of $K$ with $\sigma_+$ 
positive and $\sigma_-$ negative holes.
\item[{\rm (ii)}] The longitude $l_k\in G$ is represented by a word
$$
w_1^{-1}x_{i_1}^{\epsilon_1}w_1\ \ldots\
 w_n^{-1}x_{i_n}^{\epsilon_n}w_n,
$$
where each $w_i$ is a word in $\{x_i^{\pm 1}\}$, each $\epsilon_i=\pm 1$,
and $\sum \epsilon_i = k$.
\end{enumerate}
\end{lem}
\begin{proof} Suppose that (i) holds. By retracting the compressing disk to a 
one-dimensional spine we see that the path $l_k$, which is the boundary of 
the compressing disk, is homotopic to a product of $\sigma_+$ positive and
$\sigma_-$ negative Wirtinger generators. (ii) follows.
Conversely, we can use a homotopy between a product of $\sigma_+$ positive 
and $\sigma_-$ negative Wirtinger generators and the path $l_k$ to construct
a compressing disk with $\sigma_+$ positive and $\sigma_-$ negative holes.\end{proof}

Thus we can reinterpret the framing function $n_K\co \Z \to \N$ as follows.
For $k\in \Z$ we let 
$$
n_K(k)=\min\{n\in \N\ |\ l_k \hbox{ is represented by a word in the above 
form}\}.
$$
Note that $n_K(k)$ is independent of the choice of $\gamma$. Roughly speaking, 
we are trying to express the longitude with linking number $k$ as a shortest 
possible product of conjugates of the generators $x_1,\ldots,x_m$.
%

We call a finite presentation of a knot group in which all generators
are Wirtinger generators a {\sl Wirtinger presentation}.
The Cayley graph $\Gamma$ associated to such a presentation has
natural `layers' corresponding to the
elements' images under the natural map to $H_1(S^3\basl K)\cong \Z$.
Multiplying a given element of the knot group by any conjugate of a 
positive or negative Wirtinger generator corresponds to `stepping one layer
up' or `down' respectively in $\Gamma$. 
There are $k$ conjugates algebraically in a word representing $l_k$.
In the cases below, we shall be trying to use as few as possible,
so we want to avoid taking steps `down' in $\Gamma$, ie
using negative conjugates.  

Consider the {\sl positive cone} from the 
identity in $\Gamma$, the set of elements which may be written as a
product of positive conjugates.  We want to know how close $l_k$  
gets to this cone as $k$ increases; if we can show that, however
large $k$ gets, $l_k$ still requires steps down, 
we will obtain a negative upper bound on the natural framing of $K$.
For the left-hand trefoil, 
for example, we prove that $l_k$ requires   
at least 2 steps down 
for any $k$, corresponding to 4 extra holes of a compressing disk;
hence $n(k)\geqslant k+4$, and, since we know that for the trefoil $n(-2)=2$, 
we see $\nf(T(3,2))=-2$.

We shall show that $\nu(T(p,q))=-(p-1)(q-1)$ in two theorems, first
for the case $q=2$ and then for $q\geqslant2$.  The first proof is
really a degenerate case of the second, but we introduce the ideas
used in both theorems in the simpler context of $q=2$ and leave the
additional calculation to the second case, where it first becomes
necessary.

\begin{theorem}\sl The natural framing of $T(p,2)$ is $-(p-1)$; indeed,
its framing function is given by $n(k)=(p-1)+|k+(p-1)|$.
\end{theorem}
\begin{figure}[htb]
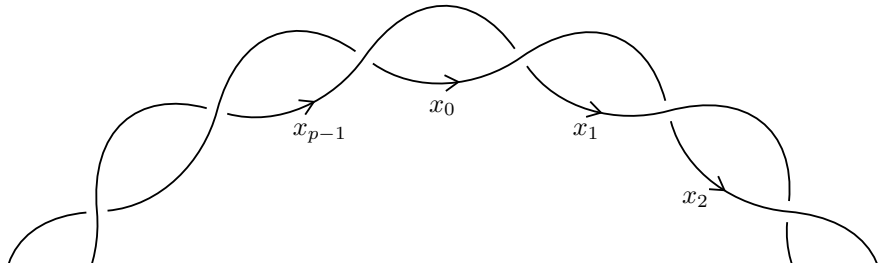

\vspace{-6mm} 
\vbox to1.8truein{%
  \vfill
  $$
  \Mon
  \mdiag{%
    .7 setlinewidth
    gsave -400 -50 moveto 400 -50 lineto 400 100 lineto -400 100 lineto
    clip newpath /ang 20 def /d -140 def /r 170 def
    /l 20 def /l2 31 def /c 4 def
    gsave 1 1 scale -4.5 1 2.5
      {90 exch ang mul sub /th exch def 0 d r th rt add2
      c th 135 sub rt add2 2 copy moveto l th 135 sub rt add2
      0 d r th ang sub rt add2 2 copy 6 2 roll 2 copy l th ang sub 135 add rt
      add2 4 2 roll curveto l2 th ang sub 45 sub rt add2 0 d r th ang 2 mul sub
      rt add2 2 copy l2 th ang 2 mul sub 45 add rt add2 4 2 roll
      c th ang 2 mul sub 45 add rt add2 curveto}
    for 0 setlinecap stroke grestore
    /i 23 def
    0 d r i sub 90 rt add2 moveto (0) show
    0 d r i sub 90 ang sub rt add2 -2 0 add2 moveto (1) show
    0 d r i sub 90 ang 2 mul sub rt add2 -5 0 add2 moveto (2) show
    0 d r i sub 90 ang add rt add2 -7 0 add2 moveto (p-1) show
    grestore
    -1 1 2 {gsave 0 d translate ang mul neg rotate 0 157.65 translate
      5 rotate 1 3 moveto 6 0 lineto 1 -3 lineto stroke grestore} for
  }%
  \Mrl{0}{\ninepoint\mlap{$x_0$}}%
  \Mrl{1}{\ninepoint$x_1$}%
  \Mrl{2}{\ninepoint$x_2$}%
  \Mrl{p-1}{\ninepoint$x_{p-1}$}%
  \Moff
  $$
  \vfill
}%
\vspace{-.8cm}
\caption{Generators of the knot group of the $(p,2)$--torus knot}\label{p2knot}
\end{figure}

\begin{proof} Fix $p$.  Draw $T(p,2)$ in the usual way in a diagram with $p$
negative crossings and $p$--fold symmetry as in Figure \ref{p2knot}.
Label the overcrossing arcs
$x_0$, $x_1$, $\ldots,$ $x_{p-1}$ clockwise around the diagram.  Then
the fundamental group of the complement is given by
$$
G={<}x_0,x_1,\ldots,x_{p-1}\>|\>x_0x_{p-1}=\ldots=x_2x_1=x_1x_0{>},
$$
and the word
$$
x_0^{k+p}x_1\inv x_3\inv\ldots\>x_{p-2}\inv x_0\inv x_2\inv\ldots\>x_{p-1}\inv
$$
represents the longitude $l_k\in G$.
Notice that $n_K(-(p-1))=p-1$; this follows from the existence of the disk
drawn in Figure \ref{mtgdisk}, which has $p-1$ negative holes and no positive 
ones, and the fact that $n(k)\geqslant|k|$.
\begin{figure}[htb]
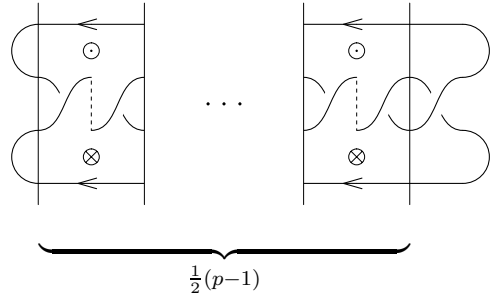

\vskip4truemm
\vbox to27truemm{%
  \Mon
  $$
  \mdiag{%
gsave 90 rotate 
gsave -3 8 tru
1 0 1 0 180 marcu 5 0 1 0 180 marcu stringline
vert -venobl vert -8 -2 tru
vert -venotr vert -8 -8 tru
vert -venobl vert -8 -2 tru
vert -venotr vert -8 -2 tru
vert -ve vert -8 -2 tru
1 0 1 180 360 marcu 5 0 1 180 360 marcu stringline
grestore
-3 6 90 arrw 3 6 90 arrw -2 6 down 2 6 up
-1 6 mvu 2 0 rlnu intersectionline
-3 -4 90 arrw 3 -4 90 arrw -2 -4 down 2 -4 up
-1 -4 mvu 2 0 rlnu intersectionline
90 0 1 3dots 
/w 3.8 def /z .3 def
gsave
  w neg 3.9 mvu w 3.9 lnu
  w neg -1.9 mvu w -1.9 lnu
  .2 unit mul setlinewidth 1 setgrey stroke
grestore
w neg 8 mvu w 8 lnu
w neg 4 mvu w 4 lnu
w neg -2 mvu w -2 lnu
w neg -6 mvu w -6 lnu symbolline
-4 1 mvu (0) show
grestore
%
  }%
  $$
  \Moff
  \vfill
  \hbox to\hsize{\hfill
    $\underbrace{\kern140pt}_{%
      {1\over2}(p-1)}$\kern20pt\hfill}%
}%
\vspace{-0.1cm}
\caption{A disk with $p-1$ negative holes and no positive ones}\label{mtgdisk}
\end{figure}

Hence $n(k)=-k$ for $k\leqslant -(p-1)$.
If we could show that, given
$k$, there was a possibly larger integer (which instead we call $k$)
such that $n(k)\geqslant k+2(p-1)$, we would know the entire framing 
function.

We aim to capture something of the geometry of the Cayley graph
$\Gamma$ of $G$ with the above presentation, and show that $l_k$ is
always `hard to get to'.  As an example, take $p=3$ and
consider just a small portion
of $\Gamma$, those vertices which may be written as the product
of at most two (positive) generators.
Three of these words coincide as group elements;
apart from that, they are all different.  This portion may be embedded
in $\R^3$ as shown in Figure \ref{cayleyg}.
\begin{figure}[htb]
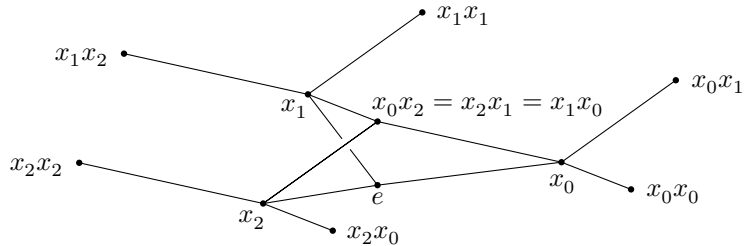

\vglue8truemm
\vbox to37truemm{%
  \vfill\Mon
  $$
  \mdiag{%
/hh 1.2 def /R 7 def /r .4 6 mul 7 div def /th -8 def
/e {0 hh neg} def
/x0 {R th rt r mul} def
/x1 {R th 120 add rt r mul} def
/x2 {R th 240 add rt r mul} def
/x0x0 {x0 R th 60 sub rt r mul add2 hh add} def
/x0x1 {x0 R th 60 add rt r mul add2 hh add} def
/x1x1 {x1 R th 60 add rt r mul add2 hh add} def
/x1x2 {x1 R th 180 add rt r mul add2 hh add} def
/x2x2 {x2 R th 180 add rt r mul add2 hh add} def
/x2x0 {x2 R th 60 sub rt r mul add2 hh add} def
/x0x2 {0 hh} def
/bl {gsave tru 0 0 1.2 0 360 marc fill grestore} def
0 -1.1 tru
e mvu x1 lnu stringline
x2 mvu x0x2 lnu gsave 5 setlinewidth 1 setgrey stroke grestore
x0x2 mvu x2 lnu e lnu x0 lnu x0x2 lnu x1 lnu
x0x0 mvu x0 lnu x0x1 lnu x1x1 mvu x1 lnu x1x2 lnu x2x2 mvu x2 lnu x2x0 lnu
stringline
e bl x0 bl x1 bl x2 bl x0x0 bl x0x1 bl x1x1 bl x1x2 bl x2x2 bl x2x0 bl x0x2 bl
1.005 dup scale  
0 -1.9 mvu (e) show
x0 0 -.9 add2 mvu (0) show
x1 -.5 -.65 add2 mvu (1) show
x2 -.4 -.7 add2 mvu (2) show
x0x0 .45 -.25 add2 mvu (00) show
x0x1 .45 -.25 add2 mvu (01) show
x1x1 .45 -.25 add2 mvu (11) show
x2x0 .45 -.25 add2 mvu (20) show
x1x2 -.45 -.25 add2 mvu (12) show
x2x2 -.45 -.25 add2 mvu (22) show
-.3 1.65 mvu (02) show
  }%
  $$
  \Mrl e{\ninepoint\mlap{$e$}}%
  \Mrl0{\ninepoint\mlap{$x_0^{}$}}%
  \Mrl1{\ninepoint\mlap{$x_1^{}$}}%
  \Mrl2{\ninepoint\mlap{$x_2^{}$}}%
  \Mrl{00}{\ninepoint$x_0^{}x_0^{}$}%
  \Mrl{01}{\ninepoint$x_0^{}x_1^{}$}%
  \Mrl{11}{\ninepoint$x_1^{}x_1^{}$}%
  \Mrl{20}{\ninepoint$x_2^{}x_0^{}$}%
  \Mrl{12}{\llap{\ninepoint$x_1^{}x_2^{}$}}%
  \Mrl{22}{\llap{\ninepoint$x_2^{}x_2^{}$}}%
  \Mrl{02}{\ninepoint$x_0^{}x_2^{}=x_2^{}x_1^{}=x_1^{}x_0^{}$}%
  \Moff\vfill
}%
\vglue-1truemm
\vspace{-5mm}
\caption{Embedding the Cayley graph of $G$ in $\R^3$}\label{cayleyg}
\end{figure}

\medbreak\noindent
We may extend this embedding in a consistent way to the whole
of $\Gamma$.  (The exact meaning of `consistent' is given
implicitly by the definition of $\theta$ below.)
Arranged thus, $\Gamma$ projects vertically down onto an infinite
$p$--valent tree, which we can think of as the Cayley graph of
$$
I_p:=*_p\Z_2
  ={<}0,1,\ldots,p-1\>|\>00=11=\ldots=(p-1)(p-1)=e{>}.
$$
(See below for a picture of this group's Cayley graph.)
To realise this projection, we want a function $\theta\co G\rightarrow I_p$.
We shall find that $\theta$ is {\it not} a homomorphism.  We
first define $\theta$ on $X%
%
%
$, the set of words in the
symbols $x_0$, $x_1$, $\ldots,$ $x_{p-1}$ and their inverses:
$$
\theta\co x_{i_0}^{\epsilon_0}x_{i_1}^{\epsilon_1}\ldots
  x_{i_{s-1}}^{\epsilon_{s-1}}\mapsto
  \prod_{j=0}^{s-1}(i_j+h_j),
$$
where each $i_j\in\{0,1,\ldots,p-1\}$,
each $\epsilon_j=\pm1$, and the $j$th `height' is
$$
h_j={\epsilon_j-1\over2}+
  \sum_{k=0}^{j-1}\epsilon_k
$$
(all addition modulo $p$).
For example, when $p=3$, the word $x_1x_2x_1x_0^{-1}x_0x_2$ maps to $12$:

\vbox{\parskip=1pt%
\bigskip\noindent
\def\ph#1{\phantom{$\nearrow$}%
  \llap{\hbox to0pt{\hss#1\hss\phantom{$\nearrow$}}}}%
\def\no{\phantom{$\nearrow$}}%
\def\ne{$\nearrow$}%
\def\se{$\searrow$}%
\newdimen\six\six=6cm\advance\six by-2.5truemm
\def\crossout{\diagplease{-9 1.5 moveto 9 5 lineto .4 setlinewidth stroke}}%
%
%
%
\noindent
\kern\six\llap{\mlap{$i_j$}\k{13}}\ph{1}\ph{2}\ph{1}\ph{0}%
  \ph{0}\ph{2}\ \ +%
\vglue2truemm
{\begingroup\par\nointerlineskip\newdimen\backup\backup=-1.7pt
\noindent\kern\six\no\no\no\no\no\ne\par\vglue\backup
\noindent\kern\six\no\no\ne\se\ne\no\par\vglue\backup
\noindent\kern\six
  \llap{\mlap{\smash{%
  \raise1.5truemm\hbox{\k{.5}$\displaystyle\sum_{k=1}^j\epsilon_k$}}}\k{13}}%
		  \no\ne\no\no\no\no\par\vglue\backup
\noindent\kern\six\ne\no\no\no\no\no\par
\endgroup}%
\vglue2truemm
\noindent\kern\six\llap{\mlap{$h_j$}\k{13}}\ph{0}\ph{1}\ph{2}\ph{2}%
  \ph{2}\ph{0}\par%
\vglue1truemm
\noindent\kern\six$\overline{\hbox{\no\no\no\no\no\no}}$\par
\vglue-3mm
\noindent\kern\six\llap{\mlap{$i_j+h_j$}\k{13}}\ph{1}\ph{0}\crossout\ph{0}%
  \ph{2}\crossout\ph{2}\ph{2}%
  ${}=12\in I_3$.
\bigskip
}%

As this example shows, changing a word in $X$ by an elementary
expansion or reduction (that is, insertion or deletion of a pair
$x_i{x_i}\!^{-1}$ or ${x_i}\!^{-1}x_i$) does not change its image under
$\theta$, since the two adjacent letters involved have the same index
and height, and so map to a repeated element in $I_p$.  Also,
$x_{i+1}x_i$ maps to the identity at any height for any $i$, so
changing our word in $X$ by a relator of $G$ leaves its image under
$\theta$ unchanged.  Hence $\theta$ is well-defined as a map from $G$
to $I_p$.

Suppose we have a word
$$
W=\yconj1\ \yconj2\ \ldots\ \yconj{n}%
$$
which represents $l_k$,   
where $n=k+2t$ and each $y_i\in\{x_0,x_1,\ldots,x_{p-1}\}$.
By adding small positive umbrellas, say by
premultiplying $W$ by a power of $x_0$, we may assume that $k$ is a
large positive multiple of $p$.  (This is for notational convenience
and to remove a special case later.  Notice that $t$ is unaffected.)
Since $W$ represents $l_k$,  
%
$$
\vbox{%
\halign{\kern-9mm\kern8.5truemm\hfill$#$&
        $#$\hfill&      $#$\hfill&\hfill$#$\hfill&\hfill$#$\hfill&
  \hfill$#$\hfill&\hfill$#$\hfill&\hfill$#$\hfill&\ \hfill$#$\hfill\ &
  \hfill$#$\hfill&\hfill$#$\hfill&\hfill$#$\hfill&\ \hfill$#$\hfill\ &
  \hfill$#$\hfill&\hfill$#$\hfill&\hfill$#$\hfill&$#$\hfill\cr
\theta(W)&{}=\theta\bigl(&(&x_0&x_0&\,\ldots\,&x_0&)^{\kp+1}&%
  x_1^{-1}&x_3^{-1}&\,\ldots\,&x_{p-2}^{-1}&x_0^{-1}&x_2^{-1}&%
  \,\ldots\,&x_{p-1}^{-1}&\bigr)\cr
\noalign{\medskip\medskip}%
&$\llap{\hbox to1em{\hss$i_j$\hss}}$&(&0&0&&0&)^{\kp+1}&%
  1&3&&p-2&0&2&&p-1&\cr                
&$\llap{\hbox to1em{\hss$h_j$\hss}}$&(&0&1&&p-1&)^{\kp+1}&%
  p-1&p-2&&\frac{p+1}2&\frac{p-1}2&\frac{p-3}2&&0&\cr
\noalign{\vglue-1mm}%
&&$\rlap{$\overline{\hbox to7cm{\hfill}}$}$&&&&&&&&&&&&&&%
  $\llap{$\overline{\hbox to7cm{\hfill}}$}$\cr
&$\llap{\hbox to1em{\hss$i_j+h_j$\hss}}$&(&0&1&\,\ldots\,&p-1&)^{\kp+1}&%
  0&1&\,\ldots\,&\frac{p-3}2&\frac{p-1}2&\frac{p+1}2&\,\ldots\,&p-1&\cr
\noalign{\medskip\medskip}%
&\rlap{${}=\bigl(0\,1\ \ldots\ (p-1)\bigr)^{\kp+2}\in I_p$.}&&&&&&&&&&&&&&&\cr
}%
}%
$$
If the Cayley graph of $I_p$ with the above generators is
drawn in the plane with edges
consistently labelled 0 to $(p-1)$ anticlockwise round each vertex, we find
that $\theta(l_k)$ turns sharp right at every step, and
follows the boundary of one of the infinite
complementary regions.  For example,
in the case $p=3$ and $k=0$ we have $\theta(l_k)=012012$,
so the graph looks as in Figure \ref{hyperb}.
\begin{figure}[htb]
\def\tweak#1#2#3#4{\Mrl{#1}{\smash{\raise#4truemm\hbox to0pt{%
  \k{#3}\mlap{\ninepoint$#2$}\hss}}}}%
\def\tweaks#1#2#3#4{\Mrl{#1}{\smash{\raise#4truemm\hbox to0pt{%
  \k{#3}\mlap{\ninepoint$\scriptstyle#2$}\hss}}}}%
\def\tweakss#1#2#3#4{\Mrl{#1}{\smash{\raise#4truemm\hbox to0pt{%
  \k{#3}\mlap{\ninepoint$\scriptscriptstyle#2$}\hss}}}}%

\vbox to3.8in{%
  \vfill\hbox to\hsize{\hfill
  \Mon
  \input ppd4tex.dummy
\includegraphics{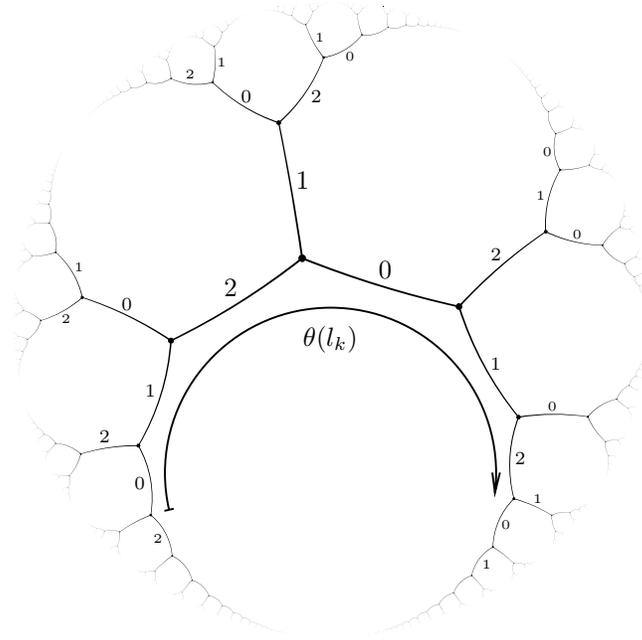}
  \special{ps: end enddiag }%
  \Mrl{theta}{\ninepoint\mlap{$\theta(l_k)$}}%
  \tweakss{28}2{-1.2}{-1}%
  \tweaks{16}0{-1.5}{-.8}%
  \tweaks92{.2}{.8}%
  \tweaks41{-1.1}{.5}%
  \tweaks80{.8}{.3}%
  \tweakss{24}2{1}{-1.3}%
  \tweakss{15}1{1.1}{-.2}%
  \tweak12{-.5}{1.2}%
  \tweak01{1.7}{-.9}%
  \tweaks20{1.1}{.1}%
  \tweakss{11}2{.2}{.7}%
  \tweakss61{.8}{-.8}%
  \tweaks52{1}{-1.9}%
  \tweakss{17}1{1.1}{-.4}%
  \tweakss{10}0{.3}{-1.6}%
  \tweak30{.4}{1.1}%
  \tweaks72{-1.4}{.1}%
  \tweakss{23}1{-1.2}{-.7}%
  \tweakss{39}0{-.9}{-.8}%
  \tweakss{14}0{-.2}{.6}%
  \tweaks{12}1{1.2}{0}%
  \tweakss{33}0{-1}{.5}%
  \tweaks{21}2{1.4}{-.3}%
  \tweakss{37}1{0}{.6}%
  \tweakss{58}0{1.1}{-.7}%
  \tweakss{95}1{1.1}{-1.1}%
  \Moff
  \hfill}\vfill
}%
\vspace{-.5cm}
\caption{The Cayley graph of the group $I_p$} \label{hyperb}
\end{figure}
%


\noindent
This suggests
defining the {\df angle} $a(v)$ of a non-trivial reduced word $v$ in
the symbols
$\{0,1,\ldots,p-1\}$, which we think of as `turn-right-ness', as
$$
a(i_0i_1\ldots i_{s-1}):=\sum_{j=1}^{s-1}(p-2d_j),
$$
$$
\hbox{ \ where }%
  i_{j-1}+d_j\equiv i_j\hbox{ (mod $p$) and }d_j\in\{1,2,\ldots,p-1\}.
$$
Thus each step $i_{j-1}$ to $i_j$ contributes between $-p+2$ and $p-2$ to
the angle, and
the more often and more sharply a word `turns right', the greater
this angle.  Define the {\sl angle} of an element of $G$ as the angle
of its reduced image under $\theta$.  Since $\theta(l_k)$ turns sharp right
$(k+2p-1)$ times, the angle is $(k+2p-1)(p-2)$; this will prove
unusually high for its exponent sum.  
Notice that the angle of an element of
$I_p$ is unchanged if we cycle the generators in its expression modulo
$p$; that is,
$$
a(i_0i_1\ldots i_{s-1})=a\Bigl((i_0+1)(i_1+1)\ldots(i_{s-1}+1)\Bigr).
$$

%

Let $W_0$ be the large power of $x_0$ we multiplied our original word by, and
define inductively $W_i=W_{i-1}\yconj i$, so $W_n=W$.
We ask how the angles of these initial segments increase with
$i$.  Suppose we are stepping from $W_{i-1}$ to $W_i$.  Since the
angle is defined only for reduced words, we must cancel
$\theta(W_i)$ down to its reduced form; we therefore pretend that
$W_i$ is a word in the letters $\{x_j^{\pm1}\}$, and that
$\theta(W_i)$ is a word in the generators of $I_p$ (the letter-by-letter image
of $W_i$), and do this reduction in stages.

First, perform all elementary reductions on $\theta(W_{i-1})$, the initial
segment.  The total angle for this segment is $a(W_{i-1})$.  Also, let
$v'$ and $v$ be the reduced forms of the words obtained from the
segments $w_i^{-1}$ and $w_i$, respectively---notice that $v$ is the
reverse of $v'$ cycled by $\epsilon_i$, so
$a(v')+a(v)=0$
(ie the angles for these segments cancel).

Next, consider the image of the whole segment $\yconj i$.
We currently have this
in the form $v'zv$, say, for some $0\leqslant z<p$.  This looks
somewhat like a conjugate of $z$ in $I_p$, except that
the generators in $v$ have been cycled by $\ei$ modulo $p$.
The words $v$ and $v'$ may already be trivial.  If not, write $z'$ for the last
letter of $v'$; then the first letter of $v$ is $(z'+\ei)$.  If
$z$ equals $z'$ or $(z'+\epsilon_i)$, we may shorten $v$ and $v'$ and
change $z$, and still have an expression of the same form (and notice
that two cancelling angles have been removed from the word).  For
example, in $I_3$, $(0212)0(0201)=(021)2(201)=(02)1(01)$.  Thus
$\theta(\yconj i)$ cancels down to one of the four forms below.  The
only contributions to the total angles for these segments come from
the steps either side of $z$, and so are as shown.
%

\nointerlineskip\goodbreak
$$
\vbox{%
\halign{\k{0}
  \hfil#&{\csc\ #}:\k7&\hfil#&%
  #\hfil\k7&\hfil#&#\hfil\k7&\hfil#&#\hfil\cr
case&a&$\ei={}$&$+1$&$\theta(\yconj i)={}$&$z$&$a(\yconj i)={}$&0\cr
\noalign{\smallskip}%
&b&&$+1$&&$v'zv$&&$-2$\cr
\noalign{\smallskip}%
&c&&$-1$&&$z$&&0\cr
\noalign{\smallskip}%
&d&&$-1$&&$v'zv$&&2\cr
}%
}%
$$
The only contribution to the angle of the whole word $\theta(W_i)$
we have not yet considered comes from the boundary between the
images under $\theta$ (reduced as described) of $W_{i-1}$ and
$\yconj i$.  If there is no cancellation at this position in the word,
the contribution to the angle of $W_i$ from this junction is at most 
$(p-2)$.  If
the second part of $\theta(W_i)$ is fully absorbed by the first, the same
is true, since each elementary reduction before the last removes two
equal and opposite contributions to the angle, and the last removes
just one, which is certainly at least $(-p+2)$.  Notice that in cases
{\csc a} and {\csc c} these are the only two possibilities.  Notice
also that we may assume that the first part is never completely
absorbed by the second, by increasing the framing beforehand as
described above.

This leaves the case where there is
partial cancellation at this boundary.  We
have so far reached
$$
\underbrace{\ldots j\,k_0k_1\ldots k_{s-1}}%
  _{\hbox{\phantom{$\scriptstyle\yconj i$}%
    $\scriptstyle\theta(W_{i-1})$\phantom{$\scriptstyle\yconj i$}}}%
\k{-1}\underbrace{k_{s-1}\ldots k_1k_0l\ldots}%
  _{\hbox{$\scriptstyle\theta(\yconj i)$}}\,,
$$
say, where $j\not=l$.  Cancelling this down, and ignoring pairs of equal and
opposite angles which vanish in the process, we lose the angle
contributions from $jk_0$
and $k_0l$, and gain instead that from $jl$.  This changes the angle by
$\pm p$.  Notice, however, that if $j+1=k_0$ the angle must decrease.

In summary, then, changes from $W_{i-1}$ to $W_i$ of types {\csc a},
{\csc b}, or {\csc c} increase the angle by at most $(p-2)$, and those
of type {\csc d} increase it by at most $(p+2)$.  The angle of $l_k$
is unusually high for its framing $k$; to see this, let ${\rm esum}(W_i)$
be the sum of the exponents in $W_i$, ie the linking number of
a path representing $W_i$ with $K$, and consider
$$
c(W_i):=a(W_i)-(p-2)\bigl({\rm esum}(W_i)-1\bigr).
$$
This is only increased in cases {\csc c} and {\csc d}, by
$(2p-4)$ and $2p$, respectively.  Recall that $a(l_k)=(k+2p-1)(p-2)$,
so $c(W_n)=2p(p-2)$, so we require at least $(p-2)$ {\csc c}s or {\csc
d}s and the same total number of {\csc a}s and {\csc b}s.  Since our word
$W$ consisted of $(k+t)$ positive conjugates and $t$ negative conjugates,
we see that $t\geqslant p-2$.

Suppose $t=p-2$.  Then $c$ must increase by a full $2p$ for
each negative conjugate, so
there are no steps of type {\csc c}.  Let $i$ be the first time we
encounter a case other than {\csc a}.  All previous steps {\it must}
have been multiplication by $x_0$ in $G$ (since the word must turn
sharp right each time to avoid decreasing $c(W_i)$).  Hence the word
so far is just a power of $x_0$.  But then at the boundary between
$W_{i-1}$ and $\yconj i$ we cannot increase $a$ by $p$, so $c$ must
decrease.  Hence $t>p-2$, so $n\geqslant k+2(p-1)$ and we are
done.\end{proof}

\begin{theorem}\sl The natural framing of $T(p,q)$ is $-(p-1)(q-1)$;
indeed, its framing function is given by
$n(k)=(p-1)(q-1)+\bigl|k+(p-1)(q-1)\bigr|$.
\end{theorem}
\begin{figure}[htb]
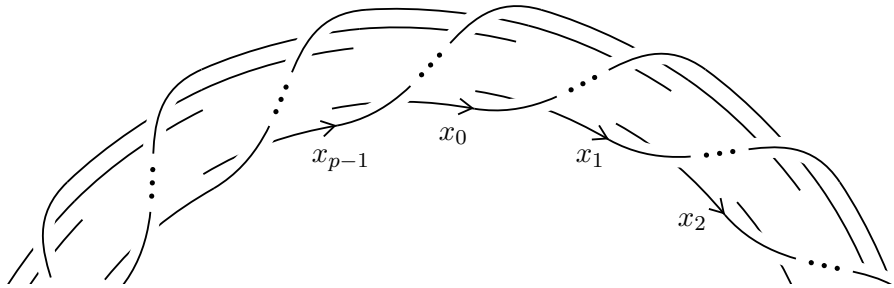

\vspace{-6mm}  
\vbox to2truein{%
  \vfill
  $$
  \Mon
  \mdiag{%
    .7 setlinewidth
    gsave -400 -54 moveto 400 -54 lineto 400 100 lineto -400 100 lineto
    clip newpath /q 6 def /ang 20 def /angstep 7 q div def
    /d -138 def /or 190 def /ir 150 def /rstep or ir sub q div def
    /l 18 def /l2 20 def /ang2 6 def /c 4 def
    gsave -1 1 scale -4.5 1 2.5
      {gsave 0 d translate ang mul neg rotate
      or rstep sub 90 angstep sub rt 2 copy moveto
        l ang2 angstep sub rt add2 or 90 ang sub
        rt 2 copy l2 180 ang sub ang2 add rt add2 4 2 roll curveto
      or rstep 2 mul sub 90 angstep 2 mul sub rt 2 copy moveto
        l ang2 angstep 2 mul sub rt add2 or rstep sub 90 ang sub angstep sub
        rt 2 copy l2 180 ang sub angstep sub ang2 add rt add2 4 2 roll curveto
      ir 90 angstep q mul sub rt 2 copy moveto
        l ang2 angstep q 1 sub mul sub rt add2 ir rstep add 90 ang sub
        angstep q mul sub rt 2 copy l2 180 ang sub angstep q mul sub
        ang2 add rt add2 4 2 roll curveto
      0 setlinecap stroke
      or 90 rt 2 copy moveto
        l 2 mul ang2 rt add2 ir 90 ang sub angstep q mul sub rt
        2 copy l2 2 mul 180 ang sub angstep q mul sub
        ang2 add rt add2 4 2 roll curveto
      gsave 1 setgrey 0 setlinecap c 2 mul setlinewidth stroke grestore
      gsave [39 22 100] 0 setdash 0 setlinecap stroke grestore
      [0 46 .01 5 .01 5 .01 100] 1 setdash .01 setlinewidth outlinestroke
      .2 setlinewidth stroke grestore}
    for -2 1 1 {gsave 0 d translate ang mul neg rotate -13 149.3 translate
      7 rotate 11 3 moveto 6 0 lineto 11 -3 lineto stroke grestore} for
    grestore
    /i 12 def
    0 d ir i sub 90 rt add2 moveto (0) show
    0 d ir i sub 90 ang sub rt add2 -2 0 add2 moveto (1) show
    0 d ir i sub 90 ang 2 mul sub rt add2 -5 0 add2 moveto (2) show
    0 d ir i sub 90 ang add rt add2 -7 0 add2 moveto (p-1) show
    grestore
  }%
  \Mrl{0}{\mlap{$x_0$}}%
  \Mrl{1}{$x_1$}%
  \Mrl{2}{$x_2$}%
  \Mrl{p-1}{$x_{p-1}$}%
  \Moff
  $$
  \vfill
}%

\vspace{-1cm}
\caption{Generators of the knot group of the $(p,q)$--torus knot}\label{mtgpq}
\end{figure}
%

\begin{proof} Fix $p$ and $q$ coprime with (say) $p>q$.
From the diagram for $T(p,q)$
shown in Figure \ref{mtgpq} we obtain the presentation
$$
G={<\,}x_0,\,x_1,\,\ldots,\,x_{p-1}\>|\>x_{q-2}\ldots x_0x_{p-1}=%
  \ \ldots\ =x_q\ldots x_2x_1=x_{q-1}\ldots x_1x_0{\,>}
$$
for $\pi_1(S^3\backslash K)$.  Now $l_k$ is the group element represented by 
the word
$$
x_0^{k+p(q-1)}x_1^{-1}x_2^{-1}\ldots x_{q-1}^{-1}x_{q+1}^{-1}\ldots
  x_{2q-1}^{-1}x_{2q+1}^{-1}\ \ \ldots\ \ x_{p-1}^{-1}.%
$$
It is less clear geometrically
in this more general case that we can find a disk corresponding to
the value $n\bigl(-(p-1)(q-1)\bigr)=(p-1)(q-1)$
of the framing function of $T(p,q)$.
Seen algebraically, however, we are trying to write
the above word
as a product of negative conjugates (again,
all indices modulo $p$) when $k=-(p-1)(q-1)$.  There are $(q-1)$ $x_0^{-1}$s in
this expression.  Conjugate each other inverse of
a generator by $x_0^{-1}$ to the
power of the number of $x_0^{-1}$s after it; then elementary reduce the $x_0$s
and $x_0^{-1}$s.  This gives an expression of the required form.  For
example, for $T(5,3)$,
%
%

\goodbreak
\vbox{%
$$
\def\,{\kern.04em}%
\def\>{\kern.08em}%
\def\ {\kern.2em}%
\halign{\k{3}\hfil#&$x_0^2\>x_1^{-1}$#\hfil&$\>x_2^{-1}$#\hfil&
  $\>x_4^{-1}$\hfil#\hfil&$x_2^{-1}$#\hfil&$x_3^{-1}\,
  x_0^{-1}\ \ x_1^{-1}\ \ x_3^{-1}\ \ x_4^{-1}$#\hfil\cr
&&&$x_0^{-1}$&&\cr
${}={}$&$\,x_0^{-2}\ \ x_0^2$&$x_0^{-2}\ \ x_0^2$&%
  $\,x_0^{-2}\ \ x_0^{}\,$&$\,x_0^{-1}\ \ x_0^{}\,$&.\cr
}%
$$
}%

We would like a version of $\theta$ for
this case; it must map a sequence of generators with
gradually decreasing index of length $q$ (rather than 2) to the
identity.  This suggests defining $\theta$ from $G$ to
$$
\pq:={<\,}0,\ 1,\ \ldots,\ p-1\>|\>0^q=1^q=\ldots=(p-1)^q=e{\,>}.
$$
As before, we first define $\theta$ on $X$, the set of words in the
symbols $x_0$, $x_1$, $\ldots,$ $x_{p-1}$ and their inverses:
$$
\theta\co x_{i_0}^{\epsilon_0}x_{i_1}^{\epsilon_1}\ldots
  x_{i_{s-1}}^{\epsilon_{s-1}}\mapsto
  \prod_{j=0}^{s-1}(i_j+h_j)^\ei,
$$
where each $i_j\in\{0,1,\ldots,p-1\}$,
each $\epsilon_j=\pm1$, and the height is again given by
$$
h_j={\epsilon_j-1\over2}+
  \sum_{k=0}^{j-1}\epsilon_k
$$
(all addition modulo $p$).
Notice that this time we keep track of inverses of 0, 1,
etc---in the involutary case, this was unnecessary.

Adjacent letters in a word in $X$ with the same index but of opposite
sign have the same height, and so map to a cancelling pair in \pq.
Applying a relator of $G$ to a word $v$ in $X$ replaces, in $\theta(v)$,
the $q$th power of some generator of \pq\ by that of another.
Therefore $\theta$ is well-defined as a map from $G$ to \pq.

Suppose we take a framing disk.
The longitude it follows is
$$
l_k=x_0^{k+p(q-1)}x_1^{-1}x_2^{-1}\ldots x_{q-1}^{-1}x_{q+1}^{-1}\ldots
  x_{2q-1}^{-1}x_{2q+1}^{-1}\ \ \ldots\ \ x_{p-1}^{-1}%
$$
(all indices modulo $p$), where $k$ may be assumed to be a large
positive multiple of $p$ by adding positive umbrellas as necessary.
This maps under $\theta$ to
$$
\begin{array}{r@{}c@{}l}
\theta(l_k)&{}={}&\bigl(01\ldots(p-1)\bigr)^{\kp+q-1}%
			0^{-(q-1)}1^{-(q-1)}\ldots(p-1)^{-(q-1)}\\
	   &{}={}&\bigl(01\ldots(p-1)\bigr)^{\kp+q}.
\end{array}
$$
The `angle' alone
is now too small to give a tight bound on the natural framing
of $T(p,q)$; we therefore introduce new ideas to make up the difference.

Since the negative powers of generators of \pq\ cancel so neatly in
$\theta(l_k)$, we want this cancelling to `score extra'.  An element
of \pq\ may be written in the {\df standard form}
$$
\matrix{%
\pretend i_0^{e_0}i_1^{e_1}\ldots i_{s-1}^{e_{s-1}},\cr
\pretend\hbox{with each $i_j$ a generator, $s$ minimal,
  and each $e_j\in\{1,2,\ldots,q-1\}$.}\cr
}%
$$
Adjacent $i_j$s are then different.  This allows us to define
{\df angle} for such a standard form of a non-trivial word,
much as before, by
$$
a(i_0^{e_0}i_1^{e_1}\!\ldots i_{s-1}^{e_{s-1}}):=
  \sum_{j=1}^{s-1}(p-2d_j),
$$
$$
\hbox{where }%
  i_{j-1}+d_j\equiv i_j\hbox{ (mod $p$) and }d_j\in\{1,2,\ldots,p-1\}.
$$
In addition, define the {\df
pseudo-exponent} of $w\in G$ to be the sum of the exponents
in the standard form of $\theta(w)$.  The {\df exponent sum of $w$}
is the sum of exponents when $w$ is written in terms of the $x_i$s.
Then the {\df excess exponent},
$e(w)$, is defined
to be the pseudo-exponent minus the exponent sum of
$w$.  Notice that $e(w)$ is a multiple of $q$.  Finally, to register how close
the excess exponent is to changing, define the {\df internal angle}
$\iota$ of a word in standard form by
$$
\iota(i_0^{e_0}i_1^{e_1}\ldots i_{s-1}^{e_{s-1}}):=
  \sum_{j=0}^{s-1}(e_j-1).
$$
Then the replacement we use for angle, which we call {\df angle\k{.1}$'$}
and denote $a'$, is given by
$$
a'(w):=a(w)\,-\,p\times\iota(w)\,+\,\frac1q(q-2)p\times e(w).
$$

Notice that $a'$ is well-defined as a map from $G$ to $\Z$,
since any element of \pq\ has a unique standard form.
Angle$'$ is {\it not} defined on \pq,
even though angle and internal angle are,
because the excess exponent of an element is undefined.
The angle$'$ of $l_k$ is 
$$
\begin{array}{r@{}c@{}l}
a'(l_k)&{}={}&(k+pq-1)(p-2)-0+\frac1q(q-2)p(k+pq-k)\\
       &{}={}&(k-1)(p-2)+2p\bigl((p-1)(q-1)-1\bigr),\\
\end{array}%
$$
and this will turn out to be large enough to give the
wanted tight upper bound on the natural framing.


Suppose we have a word
$$
W=\yconj1\ \yconj2\ \ldots\ \yconj{n}%
$$
which represents $l_k$,             
where $n=k+2t$ and each $y_i\in\{x_0,x_1,\ldots,x_{p-1}\}$.
As before, we may assume that $k$ is a
large positive multiple of $p$.  Let $W_0$ be the large power of $x_0$ which
$W$ starts with, and define $W_i=W_{i-1}\yconj i$.

Although these are really group elements, we write them out literally
in the letters $\{x_i^{\pm1}\}$ and pretend that they are just words.
We use the literal version of $\theta$ to map them to \pq.  Only then
do we do any simplification, ie elementary reduction and
multiplication by $q$th powers of generators.  In this way, we can
monitor the total angle$'$ of the result.  We shall often calculate the
angle$'$ of a segment of a word---this is again simply the angle$'$ of the
segment taken in isolation, not counting contributions to angle or
internal angle from either end.  The angle$'$ of a set of segments
of a word is the sum of the angle$'$s for each segment.

We know that simplifying the $W_{i-1}$ segment of $W_i$ gives angle$'$
$a'(W_{i-1})$.  Consider next $w_i^{-1}$ and $w_i$.  Taken as words, these
map to, say,
$$
\displaylines{%
v':=v_0^{e_0}v_1^{e_1}\ldots v_{s-1}^{e_{s-1}}\cr 
\lwd{and}v:=(v_{s-1}+\epsilon_i)^{-e_{s-1}}\ldots%
  (v_1+\epsilon_i)^{-e_1}(v_0+\epsilon_i)^{-e_0}\cr
}%
$$
(addition modulo $p$) respectively.
We would like to show that these two expressions, taken together,
contribute a total of 0 to $a'$; however,
we must convert them to standard form before we can measure this
contribution.  (By slight abuse of notation, we continue to call them
$v'$ and $v$ through the stages of this standardisation.  We are
effectively showing that $a'(g)=-a'(g^{-1})$ for any $g\in G$.)

We may assume that each $e_j$ lies in $\{0,1,\ldots,q-1\}$
(by repeatedly multiplying by \smash{$v_j^{\pm q}$} in $v'$ and
\smash{$(v_j+\epsilon_i)^{\mp q}$} in the corresponding place in $v$),
and that none is 0 (by reducing $s$).  To finish standardisation,
we must, for each $j$,
write the negative powers of $(v_j+\ei)$ as positive powers by multiplying by
$(v_j+\epsilon_i)^q$ in the correct place in $v$.
The pseudo-exponent then becomes $qs$, so $e(w_i^{-1})+e(w_i)=qs$.
For each $j$ we have $e_j$ $v_j$s in $v'$ and
$(q-e_j)$ $(v_j+\epsilon_i)$s in $v$, so the total internal
angle, $\iota(v')+\iota(v)$,
is $(q-2)s$.  We already know the total angle is zero---the
contribution from $v_jv_{j+1}$ in $v'$ cancels with that from
$(v_{j+1}+\epsilon_i)(v_j+\epsilon_i)$ in $v$---so the total angle$'$
for these two segments is, as wanted,
$$
a'(v')+a'(v)=0\,-\,p\times(q-2)s\,+\,\frac1q(q-2)p\times qs=0.
$$

Next we consider the whole second segment of $\theta(W_i)$,
namely $\theta(\yconj i)$, where the
images under $\theta$ of $w_i^{-1}$ and $w_i$ have already been
simplified as above.  Let $z^\ei$ be the image under $\theta$ of $y_i^\ei$.
Perhaps $s=0$ in the above expressions.  If so,
$\theta(\yconj i)=z^\ei$, which is already in standard form if $\ei=+1$,
and may be written in the standard form $z^{q-1}$ when $\ei=-1$ making
$e=q$ and $\iota=q-2$; both these expressions have angle$'$ 0.

Consider the case $s>0$. By cycling the whole image modulo $p$, we may
assume, for notational convenience, that $v_{s-1}$ and $(v_{s-1}+\ei)$
are 0 and 1 in some order.  The only extra
contributions to $a'$ are those involving $z^\ei$.  We watch how
the angle$'$ changes from $a'(v')+a'(v)$ to the angle$'$ of $v'zv$
after cancellation.

Suppose first that $\epsilon_i=+1$, so we find $\ldots0z1\ldots$ in
the middle of $\theta(\yconj i)$.  If $z$ is not 0 or 1, then
$\iota$ and $e$ are
unchanged and $a$ increases by $(p-2z)+(p-2(p-z+1))=-2$ (by definition,
since $z$ and $(p-z+1)$ are both in $\{1,2,\ldots,p-1\}$).  It follows that
$a'(\yconj i)=-2$.  If $z=0$ and $e_{s-1}<q-1$, or if $z=1$ and $e_{s-1}>1$,
we find $e$ unchanged, $\iota$ increased by 1, and $a$ increased by
$(p-2)$; again, $a'(\yconj i)=-2$.  Finally, if $z=0$ and $e_{s-1}=q-1$,
or if $z=1$ and $e_{s-1}=1$, we may pretend $\theta(\yconj i)$ is really a
(possibly trivial) conjugate of 1 (the first letter of
$\theta(w_i)$) or 0 (the last letter of $\theta(w_i^{-1})$) respectively 
with smaller $s$.  This simplification removes two equal and opposite
angles, reduces $e$ by $q$, and reduces $\iota$ by $(q-2)$, leaving
$a'$ unchanged.  By induction, cancellation
of this kind terminates in one of the cases already considered.

Suppose instead that $\epsilon_i=-1$, so we find $\ldots1z^{-1}0\ldots$ in
the middle of the word
$\theta(\yconj i)$.  We may write the central $z^{-1}$ as
$z^{q-1}$, since this increases $\iota$ by $(q-2)$ and $e$ by $q$, leaving
$a'$ unchanged.  Now, if $z$ is not 0 or 1, then $\iota$ and $e$ are
unaltered and the angle $a$ increases by $(p-2(z-1))+(p-2(p-z))=2$ (as above,
since now $(z-1)$ and $(p-z)$ are in $\{1,2,\ldots,p-1\}$).  In this case
$a'(\yconj i)=2$.  If $z=1$ and $e_{s-1}>1$ (or if $z=0$ and $e_{s-1}<q-1$)
then
a block of $q$ 1s (or 0s) may be removed from the centre, leaving
$\ldots1^{e_{s-1}-1}0^{q-e_{s-1}}\ldots$
(or $\ldots1^{e_{s-1}}0^{q-e_{s-1}-1}\ldots$).
In this middle segment, $e$ is reduced by $q$, $\iota$ is reduced from
$(q-2)+(q-2)$ to $(q-3)$, and $a$ is increased by
$(p-2)$; again, $a'(\yconj i)=2$.  Finally,
if $z=1$ and $e_{s-1}=1$, or if $z=0$ and $e_{s-1}=q-1$, we may pretend
that $\theta(\yconj i)$ is really a
(possibly trivial) conjugate of $0^{q-1}$ (the first $(q-1)$ letters of
$\theta(w_i)$) or $1^{q-1}$ (the last $(q-1)$ letters of $\theta(w_i^{-1})$)
respectively  with smaller $s$.  Again, the simplification removes
two equal and opposite angles, reduces $e$ by $q$, and reduces
$\iota$ by $(q-2)$, and therefore leaves $a'$ unchanged;
hence, by induction, cancellation
of this kind terminates in one of the cases already considered.

Thus we have essentially the same set of cases as before:
%
$$
\vbox{%
\halign{\k{0}
  \hfil#&{\csc\ #}:\k7&\hfil#&%
  #\hfil\k7&\hfil#&#\hfil\k7&\hfil#&#\hfil\cr
case&a&$\ei={}$&$+1$&$\theta(\yconj i)={}$&$z$&$a'(\yconj i)={}$&0\cr
\noalign{\smallskip}%
&b&&$+1$&&$v'zv$&&$-2$\cr
\noalign{\smallskip}%
&c&&$-1$&&$z^{q-1}$&&0\cr
\noalign{\smallskip}%
&d&&$-1$&&$v'z^{q-1}v$&&2\cr
}%
}%
$$

Finally, there is the junction between $\theta(W_i)$ (say ending in $x$) and
$\theta(\yconj i)$ (say beginning with $y$) to consider,
where each part here has been written
in standard form as described above.  The total angle$'$
(ie the sum of the angle$'$s) of these segments,
taken separately, is $a'(W_i)$ or $a'(W_i)\pm2$, depending on the form of
$\theta(\yconj i)$.  If $x\not=y$,
the juxtaposition of these strings is already
in standard form; then $\iota$ and $e$ remain the same and $a$ increases
by at most $(p-2)$, so $a'$ increases by at most $(p-2)$.

Suppose, then, that $x=y$.
Some cancelling may occur at the boundary---if it does,
we continue to call the shorter (standard form) words whose product is
$\theta(W_{i+1})$ `the first word' and `the second word', even though they
change.  Say the first word ends in $\ldots x'x^j$,
where $x'\not=x$.  (We may always
assume that we can find such an $x'$, by having premultiplied $W$ by a
sufficiently high power of $x_0$ before we started.) The second may begin
with $x^kx''$, where $x\not=x''$, or may simply be of the form $x^k$.
Hence we have
$$
\hbox{either\k4}\ldots x'x^j\ x^kx''\!\ldots\k4\hbox{or}\k4\ldots x'x^j\ x^k
$$
at the boundary.  Since the words are in standard form, $j$ and $k$
are both in $\{1,2,\ldots,q-1\}$.  Compare $(j+k)$ with $q$:

\newdimen\nicestartskip
\setbox0\hbox{$j+k=q$}%
\nicestartskip=\wd0
\def\nicestart#1{\leavevmode\hbox to\nicestartskip{$j+k$\hss$#1$\hss$q$}:\k5}%

{%
\medskip
\parindent0em\hangindent2em
\nicestart<%
The juxtaposition of these strings is $\k1\ldots x'x^{j+k}x''\!\ldots,$
which is already in standard form,
so $\iota$ increases by 1 and $a'$ drops by $p$.

\medskip
\parindent0em\hangindent2em
\nicestart>%
The juxtaposed strings cancel to $\k1\ldots x'x^{j+k-q}x''\!\ldots,$
so $e$ decreases by $q$, but $\iota$ drops by $(q-1)$,
so $a'$ increases by $p$.

\medskip
\parindent0em\hangindent2em
\nicestart=%
$x^j$ and $x^k$ cancel, reducing $e$ by 1 and $\iota$ by
$(q-2)$, and hence leaving $a'$ unchanged.
Then one of the following holds.

\medskip

\parindent4em\hangindent4em
\leavevmode\llap{(i)\k3}%
There is no $x''$, so the word becomes $\k1\ldots x'\!$, $a$ drops by
at least $(-p+2)$, and $e$ and $\iota$ are unchanged,
so $a'$ increases by at most $(p-2)$.

\parindent4em\hangindent4em
\leavevmode\llap{(ii)\k3}%
$x'\not=x''$, so the string becomes $\k1\ldots x'x''\!\ldots,$
$e$ and $\iota$ are unchanged, and
$a$ (and hence $a'$) changes by either $+p$ or $-p$
(certainly $-p$ if $x'+1=x$).

\parindent4em\hangindent4em
\leavevmode\llap{(iii)\k3}%
$x'=x''$, and we can induct with shorter words.

}%

\medskip
The largest possible change to $a'$ from this boundary is therefore $p$,
as before, so we find that $a'(W_{i+1})-a'(W_i)$ is at most
$(p-2)$, $(p-2)$, $p$, and $(p+2)$ in the cases
{\csc a} to {\csc d}, respectively.
But $a'(l_k)=(k-1)(p-2)+2p\bigl((p-1)(q-1)-1\bigr)$, so if $t$ is the number
of negative conjugates in $W$ we have
$$
\displaylines{%
(k-1)(p-2)+t(p-2)+t(p+2)\geqslant(k-1)(p-2)+2p\bigl((p-1)(q-1)-1\bigr)\cr
\lwd{$\Rightarrow$}t\geqslant(p-1)(q-1)-1.\cr
}%
$$

Suppose $t=(p-1)(q-1)-1$.  Then there are no steps of type {\csc c}, since
each negative conjugate must increase $a'$ by the full $(p+2)$.
Any initial steps of type {\csc a} must be
multiplication of the word in $G$ by $x_0$, since $a'$ must increase by
$(p-2)$ each step and the only way to do this is to increase $a$ by the
largest possible amount.  But then the first step of type {\csc b} or {\csc d}
must {\it reduce} $a'$ by $p$ at the boundary between $W_{i-1}$ and
$\yconj i$ (see the case labelled (ii) above),
which is a contradiction.

Hence $t\geqslant(p-1)(q-1)$, and we are done.\end{proof}

The theorem yields many other knots for which the absolute values
of the natural framing number and the signature are
different. For instance, if $K$ is the $(3,7)$ torus knot then
$\nu(K)=-12$ and, according to \cite{Litherland}, $\sigma(K)=8$.
In this case, $|\nu(K)|>|\sigma(K)|$.  As M~Lackenby \cite{Lackenb}
has pointed out, one can use this fact to construct non-prime knots with
$|\nu(K)|<|\sigma(K)|$. For instance, for $K=T(3,7)\#T(2,-13)$ we have
$\nu(K)=12-12=0$ and $\sigma(K)=-8+12=4$.

If our conjecture on the natural framing of ${4ml+1\over 2l}$--two
bridge knots is true, then we 
obtain many more knots $K$ with $|\nu(K)|\neq|\sigma(K)|$.
For instance, for $m=l=2$ we would have that the knot $7_4$ has natural
framing $-4$, whereas its signature is $2$ (see the table in \cite{BuZi}).
Still, it is reasonable to expect that the natural framing number has some
properties similar to the signature:

{\bf Questions}\stdspace(1)\stdspace Is $\nu(K)$ even if $K$ is atoroidal, or if $K$ is a 
two-bridge knot?

(2)\stdspace Is the natural framing number of slice knots always zero? Is the natural
framing number a concordance-invariant?

(3)\stdspace Is $2 u(K)\geqslant |\nu(K)|$, where $u$ is the unknotting number?
Is $|\nu(K)|\geqslant|\sigma(K)|$ if $K$ is a prime knot? 
Is $|\nu(K)|\leqslant 2g(K)$, where $g(K)$ is the four-ball genus?
We conjecture that the answer to all three questions is No.


(4)\stdspace Is there a finite algorithm to compute the natural
framing number of a given knot?  What about two-bridge knots?


\section{A table}

The following table contains all we know about the natural framing
numbers of prime knots with up to seven crossings. For each of these
knots it states the range within which the natural framing could
possibly lie, and the value which we conjecture. For the convenience of
the reader, we have added columns for invariants which may be related to
the natural framing: the signature (taken from \cite{BuZi}),
the blackboard-framing of an alternating diagram (see \cite{Murasugi}), 
and the unknotting number (taken from \cite{KirbyPr}).

$$
\begin{array}{cccccc}
{\rm knot} & \rlap{\kern-3mm{\rm natural framing}} & & {\rm signature} & 
{\rm alternating} & {\rm unknotting}\\
 & {\rm range} & {\rm conjecture} &  & {\rm diagram} & {\rm number}\\
3_1 & 2 &  & -2 & 3 & 1\\
4_1 & 0 &  &  0 & 0 & 1\\
5_1 & 4 &  & -4 & 5 & 2\\
5_2 & [0,2] & 2 & -2 & 5 & 1\\
6_1 & 0 &  &  0 & 2 & 1\\
6_2 & [0,2] & 2 & -2 & 2 & 1\\
6_3 & 0 &  &  0 & 0 & 1\\
7_1 & 6 &  & -6 & 7 & 3\\
7_2 & [0,2] & 2 & -2 & 7 & 1\\
7_3 & [0,4] & 4 & -4 & 7 & 2\\
7_4 & [0,4] & 4 & -2 & 7 & 2\\
7_5 & [0,4] & 4 & -4 & 7 & 2\\
7_6 & [0,2] & 2 & -2 & 3 & 1\\
7_7 & 0 &  &  0 & 1 & 1\\
\end{array}
$$

\begin{proof} $3_1$, $5_1$ and $7_1$ are torus knots. $4_1$ and $6_1$ belong to the
family of twist knots with natural framing $0$. $6_3$ is amphichiral. It is 
easy to find
compressing disks of the knots $5_2$, $7_2$, $7_3$, $7_4$, $7_5$ and
$7_6$ with no ribbons, no negative clasps, only the appropriate number of 
positive clasps. (Note that the values for $5_2$, $7_2$ and $7_4$ have
been conjectured in section 1.) The only cases where the compressing disks
are not easy to imagine are $6_2$ and $7_7$.
\begin{figure}[htb]
$$
\Mon
\epsfbox{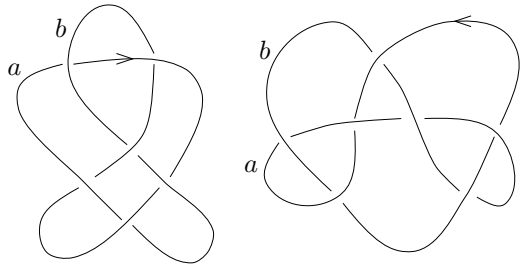}%
\Mrl a{\ninepoint$\!\!a$}%
\Mrl{a'}{\ninepoint$\!\!a$}%
\Mrl b{\ninepoint$\!b$}%
\Mrl{b'}{\ninepoint$\!b$}%
\Moff
$$
\vspace{-.8cm}
\caption{The knots $6_2$ and $7_7$}
\end{figure}

We see that $6_2$ has longitude
$$l=a^*\ bABaBAbABa\ b^*\ aBAbABaBAb
$$
(where we write $A$ and $B$ for $a\inv$ and $b\inv$, respectively);
we can write this as
\def\oo#1#2{#2}%
$$
\oo{aa}{a^2}\ bA{\underline B}aBAbA {\underline B}a\ \oo{BB}{b^{-2}}\
a{\underline B}AbABa{\underline B}Ab.
$$
If we leave out the four underlined letters we obtain the trivial word.
This proves that $n(4)=4$.
On the other hand, we can write the longitude as
$$
\oo{aaa}{a^3}\ {\underline b}ABa{\underline B}AbABa\ \oo{b}{b^1}\
{\underline a}BAbABa{\underline B}Ab,
$$
so $n(0)\leqslant 4$.  It follows that
$0\leqslant\nu(6_2)\leqslant2$.

The knot $7_7$ has longitude
$$
l=a^*\ bABaBAbaBabAbaBAbABa\ b^*\ aBAbABabAbaBabABaBAb.
$$
We can write this as
$$
\oo{aaaa}{a^4}\ bABaBAb{\underline a}BabAbaBA{\underline b}AB
{\underline a}\ \oo{bb}{b^2}\
{\underline a}BA{\underline b}ABabAbaB{\underline a}bABaB
Ab.
$$
Leaving out the six underlined letters yields the trivial word, which
proves that $n(-6)=6$.  Also, we can write the longitude as
$$
a^{-6}\ b{\underline A}BaB{\underline A}ba{\underline B}
ab{\underline A}ba{\underline B}AbA{\underline B}a\ 
b^{-6}%
\ a{\underline B}AbA{\underline B}ab{\underline A}ba
{\underline B}ab{\underline A}BaB{\underline A}b,
$$

so $n(12)=12$.  It follows that $\nu(7_7)=0$.\end{proof}

{\bf Remark}\stdspace These compressing disks were found using box-diagrams (see
\cite{nat3}). \nl

{\bf Acknowledgements}\stdspace The authors thank their respective
PhD advisors Brian Sanderson and Colin Rourke for their help and enthusiasm.
M.T.G.\ was sponsored by EPSRC, B.W.\ by a University of Warwick Graduate 
Award.

\end{document}